\documentclass[10pt]{article}

\RequirePackage{amsmath}
\RequirePackage{amsfonts}
\RequirePackage{amssymb}
\RequirePackage{hyperref}
\RequirePackage{multirow}
\RequirePackage{algorithm}
\RequirePackage{algpseudocode}
\RequirePackage{graphicx}
\RequirePackage{xcolor}
\RequirePackage{amsthm}
\RequirePackage{amsmath}
\RequirePackage{amsfonts}
\RequirePackage{amssymb}
\RequirePackage{hyperref}
\RequirePackage{multirow}
\RequirePackage{algorithm}
\RequirePackage{algpseudocode}
\RequirePackage{graphicx}
\RequirePackage{tcolorbox}
\RequirePackage{xcolor}
\allowdisplaybreaks

\usepackage{dpr_fec,ifthen,dsfont}

\usepackage{isetechreport}
\def\reportyear{21T}
\def\reportno{015}
\def\revisionno{0}
\def\originaldate{\today}
\coraltrue
\cvcrfalse

\begin{document}

\title {Inexact Sequential Quadratic Optimization for Minimizing a Stochastic Objective Function Subject to Deterministic Nonlinear Equality Constraints}

\author{Frank E.~Curtis\thanks{E-mail: frank.e.curtis@lehigh.edu}}
\author{Daniel P.~Robinson\thanks{E-mail: daniel.p.robinson@lehigh.edu}}
\author{Baoyu Zhou\thanks{E-mail: baoyu.zhou@lehigh.edu}}
\affil{Department of Industrial and Systems Engineering, Lehigh University}

\titlepage

\maketitle

\begin{abstract}
  An algorithm is proposed, analyzed, and tested experimentally for solving stochastic optimization problems in which the decision variables are constrained to satisfy equations defined by deterministic, smooth, and nonlinear functions.  It is assumed that constraint function and derivative values can be computed, but that only stochastic approximations are available for the objective function and its derivatives.  The algorithm is of the sequential quadratic optimization variety.  A distinguishing feature of the algorithm is that it allows inexact subproblem solutions to be employed, which is particularly useful in large-scale settings when the matrices defining the subproblems are too large to form and/or factorize.  Conditions are imposed on the inexact subproblem solutions that account for the fact that only stochastic objective gradient estimates are available.  Convergence results in expectation are established for the method.  Numerical experiments show that it outperforms an alternative algorithm that employs highly accurate subproblem solutions in every iteration.
\end{abstract}

\newcommand{\dTrue}{d^{\rm true}}
\newcommand{\uTrue}{u^{\rm true}} 
\newcommand{\vTrue}{v^{\rm true}} 
\newcommand{\deltaTrue}{\delta^{\rm true}}
\newcommand{\tauTrue}{\tau^{\rm true}}
\newcommand{\tautrial}{\tau^{\rm trial}}
\newcommand{\tauTruetrial}{\tau^{\rm true,trial}}
\newcommand{\alphainit}{\alpha^{\rm init}}
\newcommand{\alphasuff}{\alpha^{\rm suff}}
\newcommand{\alphamin}{\alpha^{\rm min}}
\newcommand{\alphamax}{\alpha^{\rm max}}
\newcommand{\alphavarphi}{\alpha^{\varphi}}
\newcommand{\alphatrial}{\alpha^{\rm trial}}
\newcommand{\xitrial}{\xi^{\rm trial}}
\newcommand{\Tone}{\Kcal_1}
\newcommand{\Ttwo}{\Kcal_2}
\newcommand{\Tthree}{\Kcal_3}
\newcommand{\lambdamax}{\lambda^{\rm max}}
\newcommand{\lambdamin}{\lambda_{\rm min}}
\newcommand{\low}{{\rm low}}
\newcommand{\zero}{{\rm zero}}
\newcommand{\SISQO}{{\tt SISQO}}
\newcommand{\SISQOexact}{{\tt SISQO\_exact}}
\newcommand{\xinexact}{x_{{\tt SISQO}}}
\newcommand{\xexact}{x_{{\tt SISQO\_exact}}}
\newcommand{\yinexact}{y_{{\tt SISQO}}}
\newcommand{\yexact}{y_{{\tt SISQO\_exact}}}

\section{Introduction}\label{sec.introduction}

In this paper, we consider the design, analysis, and implementation of a stochastic inexact sequential quadratic optimization (SISQO) algorithm for minimizing a stochastic objective function subject to deterministic equality constraints.  Specifically, we consider problems that may be written in the form
\bequation\label{prob.opt}
  \min_{x\in\R{n}}\ f(x)\ \st\ c(x) = 0,\ \text{with}\ f(x) = \E[F(x,\omega)],
\eequation
where $f : \R{n}\to\R{}$ and $c : \R{n} \to \R{m}$ are continuously differentiable, $\omega$~is a random variable with probability space $(\Omega,\Fcal,P)$, $F : \R{n}\times\Omega\to\R{}$, and $\E[\cdot]$ represents expectation taken with respect to the distribution of $\omega$.  Problems of this type arise in numerous important application areas.  A partial list is the following: (i) learning a deep convolutional neural network for image recognition that imposes properties (e.g., smoothness) of the systems of partial differential equations (PDEs) that the convolutional layers are meant to interpret \cite{RuthHabe20}; (ii) multiple deep learning problems including physics-constrained deep learning for high-dimensional surrogate modeling and uncertainty quantification without labeled data \cite{ZhuZabaKoutPerd19}, natural language processing with constraints on output labels \cite{NandPathAbhiSing19}, image classification, detection, and localization \cite{RaviDinhLokhSing19}, deep reinforcement learning \cite{AchiHeldTamaAbbe17}, deep network compression \cite{ChenTungVeduMori18}, and manifold regularized deep learning \cite{KumaSoumMhamHara18,TomaRose14}; (iii) accelerating the solution of PDE-constrained inverse problems by using a reduced-order model in place of a full-order model, coupled with techniques to learn the discrepancy between the reduced-order and full-order models \cite{SherRaguMoreAdamThan19}; (iv) multi-stage modeling \cite{ShapDentRusz14}; (v) portfolio selection \cite{ShapDentRusz14}; (vi) optimal power flow \cite{SummWarrMoraLyge15,VrakMathAnde14,WoodWollSheb13}; and (vii) statistical problems such as maximum likelihood estimation with constraints \cite{ChatChenMaasCarr16,Geye91}.  (For an overview of the promises and limitations of imposing hard constraints during deep neural network training, see \cite{MarqSalzFua17}.)

Popular algorithmic approaches for solving problems of the form~\eqref{prob.opt} when the objective function~$f$ is \emph{deterministic} include penalty methods~\cite{Cour43,Flet13} and sequential quadratic optimization (SQO) methods~\cite{Powe78a,Wils63}.  Penalty methods (which include popular strategies such as the augmented Lagrangian method and its variants) handle the constraints indirectly by adding a measure of constraint violation to the objective function, perhaps aided by information related to Lagrange multiplier estimates.  The resulting unconstrained optimization problem, which can be nonsmooth depending on the choice of the constraint violation measure, may be solved using a host of methods such as line search, trust region, cubic regularization, subgradient~\cite{Shor12}, or proximal methods~\cite{KaplTich98,Rock76} (with the appropriateness of the method depending on whether the constraint violation measure is smooth or nonsmooth).  It is often the case that a sequence of such unconstrained problems needs to be solved to obtain appropriate Lagrange multiplier estimates and/or to identify an adequate weighting between the original objective $f$ and the measure of constraint violation so that the original constrained problem can be solved to reasonable accuracy.

SQO methods, on the other hand, handle the constraints directly by employing local derivative-based approximations of the nonlinear constraints in explicit affine constraints in the subproblems employed to compute search directions.  For example, so-called line search SQO methods are considered state-of-the-art for solving deterministic equality constrained optimization problems~\cite{Han77,HanMang79,Powe78b}.  During each iteration of such a line search SQO method, a symmetric indefinite linear system of equations is solved, followed by a line search on an appropriate merit function to compute the next iterate.  (Here, the linear system can be seen as being derived from applying Newton's method to the stationarity conditions for the nonlinear problem, and for this reason in the setting of equality constrained optimization, SQO methods are often referred to as Newton methods.)  For large-scale problems, factorizing the matrix in this linear system may be prohibitively expensive, in which case it may be preferable instead to apply an iterative linear system solver, such as MINRES~\cite{PaigSaun75}, to the linear system.  This, in turn, opens the door to employing inexact subproblem solutions that may offer a better balance between per-iteration and overall computational costs of the algorithm for solving the original nonlinear problem.  Identifying appropriate inexactness conditions that ensure that each search direction is sufficiently accurate so that the SQO algorithm is well posed and converges to a solution (under reasonable assumptions) is a challenging task with few solutions~\cite{BiroGhat03,ByrdCurtNoce08,ByrdCurtNoce10,HeinRidz08b,HeinVice02}.

The success of SQO methods in the deterministic setting motivates us to study their extensions to the \emph{stochastic} setting, which is a very challenging task. We are only aware of three papers, namely~\cite{BeraCurtONeiRobi21,BeraCurtRobiZhou21,NaAnitKola21}, that present algorithms for solving stochastic optimization problems with deterministic nonlinear equality constraints that offer convergence guarantees with respect to solving the constrained problem (rather than, say, merely a minimizer of a penalty function derived from the constrained problem). The algorithm in~\cite{NaAnitKola21} is a line search method that uses a differentiable exact augmented Lagrangian function as its merit function, whereas   \cite{BeraCurtRobiZhou21} (resp.,~\cite{BeraCurtONeiRobi21}) is an SQO method that uses an $\ell_1$-norm (resp.,~$\ell_2$-norm) penalty function as its merit function.  All of these methods must factorize a matrix during each iteration, which may not be tractable for large-scale problems.  This motivates the work in this paper, which extends the methods in~\cite{BeraCurtONeiRobi21,BeraCurtRobiZhou21} to allow for inexact subproblem solutions, thereby making our approach applicable for solving problem~\eqref{prob.opt} in large-scale settings.

\subsection{Contributions}

The contributions of this paper pertain to a new algorithm for solving problem~\eqref{prob.opt}, which we now summarize. (i) We design a SISQO method for solving the stochastic optimization problem~\eqref{prob.opt} that is built upon a set of conditions that determine what constitutes an acceptable inexact subproblem solution along with an adaptive step size selection policy.  The algorithm employs an $\ell_2$-norm merit function, the parameter of which is updated dynamically by a procedure that has been designed with considerable care, since it is this parameter that balances the emphasis between the objective function and the constraint violation in the optimization process.  (ii) Under mild assumptions that include good behavior of the adaptive merit parameter (which can be justified as explained in the paper), we prove convergence in expectation of our algorithm.  (iii) We present numerical results that compare our SISQO algorithm to a stochastic exact SQO algorithm.  These experiments show that our SISQO algorithm benefits from our proposed inexactness strategy.

\subsection{Notation}

Let $\R{}$ denote the set of real numbers, $\R{}_{\geq p}$ (resp., $\R{}_{>p}$) denote the set of real numbers greater than or equal to (resp., strictly greater than) $p\in\R{}$, and $\N{} := \{0,1,2,\ldots\}$ denote the set of natural numbers. Let $\R{n}$ denote the set of $n$-dimensional real vectors, $\R{m\times n}$ denote the set of $m$-by-$n$-dimensional real matrices, and $\Smbb^n$ denote the set of $n$-by-$n$-dimensional symmetric real matrices. For any $p \in \N{} \setminus \{0\}$, let $[p] := \{1,\ldots,p\}$.  The $\ell_2$-norm is written simply as $\|\cdot\|$.

Our algorithm generates a sequence of iterates $\{x_k\}$ where $x_k\in\R{n}$ for all $k\in\N{}$.  For all $k \in \N{}$, we append the subscript $k$ to other quantities defined in the $k$th iteration of the algorithm, and for brevity we define $\nabla f_k := \nabla f(x_k)$, $c_k := c(x_k)$, and $J_k = \nabla c(x_k)^T$.  We refer to the range space of $J_k^T$ as $\Range(J_k^T)$ and refer to the null space of $J_k$ as $\Null(J_k)$, and recall that the Fundamental Theorem of Linear Algebra provides that these spaces are orthogonal and $\Range(J_k^T) + \Null(J_k) = \R{n}$.  Finally, recall (see, e.g., \cite{NoceWrig06}) that a \emph{primal} point $x \in \R{n}$ and \emph{dual} point $y \in \R{m}$ constitute a first-order stationary point for problem~\eqref{prob.opt} if and only if
\bequation\label{eq.stationary}
  c(x) = 0\ \ \text{and} \ \ \nabla f(x) + \nabla c(x)y = 0.
\eequation
These conditions are necessary for $x$ to be a local minimizer when the constraint functions satisfy a constraint qualification, as is assumed in the paper; see Assumption~\ref{ass.main}.

\subsection{Organization}

Our algorithm is presented in Section~\ref{sec.stochastic}.  Our convergence analysis for the algorithm is presented in Section~\ref{sec.analysis}.  The results of numerical experiments are presented in Section~\ref{sec.numerical} and concluding remarks are presented in Section~\ref{sec.conclusion}.

\section{SISQO Algorithm}\label{sec.stochastic}

Our proposed algorithm generates a sequence
\bequationNN
  \{(x_k,y_k,v_k,u_k,d_k,\delta_k,\tau_k,\alpha_k)\} \subset \R{n}\times\R{m}\times\R{n}\times\R{n}\times\R{n}\times\R{m}\times\R{}_{>0}\times\R{}_{>0},
\eequationNN
where, for all $k \in \N{}$, $(x_k,y_k)$ is a primal-dual iterate pair, $v_k$ is a \emph{normal} direction that aims to reduced infeasibility by reducing a local derivative-based model of the $\ell_2$-norm constraint violation measure, $u_k$ is a \emph{tangential} direction that aims to maintain the reduction in linearized infeasibility achieved by the normal direction while also aiming to reduce the objective function by reducing a stochastic-gradient-based quadratic approximation of the objective, $d_k := v_k + u_k$ is a full primal search direction, $\delta_k$ is a dual search direction, $\tau_k$ is a merit function parameter, and $\alpha_k$ is a step size that aims to produce $x_{k+1} \gets x_k + \alpha_k d_k$ yielding sufficient reduction in the $\ell_2$-norm merit function.  (The algorithm also generates sequences of adaptive auxiliary parameters that are introduced throughout our algorithm description.)  In the remainder of this section, we discuss each of these quantities in further detail toward our complete algorithm statement, which is provided as Algorithm~\ref{alg:stochastic_sqp_adaptive} on page~\pageref{alg:stochastic_sqp_adaptive}.

For the remainder of the paper, we make the following assumption.   

\bassumption\label{ass.main}
  Let $\Xcal\subseteq\R{n}$ be an open convex set containing the iterate sequence $\{x_k\}$ generated by any run of our algorithm.  The objective function $f:\R{n}\to\R{}$ is continuously differentiable and bounded below over $\Xcal$ and its gradient function $\nabla f:\R{n}\to\R{n}$ is Lipschitz continuous with constant $L \in \R{}_{>0}$ $($with respect to the $\ell_2$-norm$)$ and bounded over $\Xcal$.  The constraint function $c:\R{n}\to\R{m}$ $($with $m\leq n$$)$ is continuously differentiable and bounded over $\Xcal$ and its Jacobian function $J:\R{n}\to\R{m\times n}$ is Lipschitz continuous with constant $\Gamma \in \R{}_{>0}$ $($with respect to the induced $\ell_2$-norm$)$ and bounded over $\Xcal$.  In addition, for all $x \in \Xcal$, the Jacobian $J(x)$ has singular values that are bounded uniformly below by a positive real number.
\eassumption

Such an assumption is standard in the literature on deterministically constrained optimization.  Observe that it does not include an assumption that $\Xcal$ is bounded.

\subsection{Merit function}\label{sec.Merit}

Motivated by the success of numerous line search SQO methods for solving deterministic equality constrained optimization problems, our algorithm employs an exact penalty function as a merit function; in particular, it employs the $\ell_2$-norm merit function $\phi:\R{n}\times\R{}_{>0}\to\R{}$ defined by
\bequation\label{eq.merit}
  \phi(x,\tau) = \tau f(x) + \|c(x)\|,
\eequation
where $\tau$ is a positive \emph{merit parameter} that is updated adaptively by the algorithm.  (The choice of the $\ell_2$-norm in $\phi$ is not essential for our method.  Another norm could be used instead.  The choice of the $\ell_2$-norm merely makes certain calculations simpler for our presentation and analysis.)  A model $l:\R{n}\times\R{}_{>0}\times\R{n}\times\R{n}\to\R{}$ of the merit function based on $g \approx \nabla f(x)$ and $\nabla c(x)$ is given by
\bequationNN
  l(x,\tau,g,d) = \tau(f(x) + g^Td) + \|c(x) + \nabla c(x)^Td\|,
\eequationNN
with which we define the model reduction function $\Delta l:\R{n}\times\R{}_{>0}\times\R{n}\times\R{n}\to\R{}$ by
\bequation\label{eq.model_reduction}
  \baligned
    \Delta l(x,\tau,g,d)
      =&\ l(x,\tau,g,0) - l(x,\tau,g,d) \\
     :=&\ -\tau g^Td + \|c(x)\| - \|c(x) + \nabla c(x)^Td\|.
  \ealigned
\eequation
The merit function, and in particular the model reduction function \eqref{eq.model_reduction}, play critical roles in our inexactness conditions for defining acceptable search directions and in our step size selection scheme, as can be seen in the following subsections.

\subsection{Computing a search direction}\label{sec:direction}

During the $k$th iteration, the algorithm computes a normal direction $v_k\in \Range(J_k^T)$ based on the subproblem
\bequation\label{prob.normal}
  \min_{v\in\Range(J_k^T)}\ \tfrac{1}{2}\|c_k + J_kv\|^2.
\eequation
Instead of solving \eqref{prob.normal} exactly, the algorithm allows for an inexact solution to be employed by only requiring the computation of $v_k \in \Range(J_k^T)$ satisfying
\bequation\label{eq.Cauchy_decrease}
  \|c_k\| - \|c_k + J_kv_k\| \geq \epsilon_c(\|c_k\| - \|c_k + \alpha_k^c J_k v_k^c\|)
\eequation
(commonly known as the Cauchy decrease condition), where $\epsilon_c\in (0,1]$ is a user-defined constant.  In \eqref{eq.Cauchy_decrease}, $v_k^c := -J_k^Tc_k$ is the negative gradient direction for the objective of \eqref{prob.normal} at $v=0$ and $\alpha_k^c$ is the step size along $v_k^c$ that minimizes $\|c_k + \alpha J_k v_k^c\|$ over $\alpha \in \R{}$.  If $\|c_k\| \neq 0$, then under Assumption~\ref{ass.main} it follows that $\|J_k^Tc_k\| \neq 0$,
\bequation\label{eq.alphac}
  \alpha_k^c = \|J_k^Tc_k\|^2/ \|J_kJ_k^Tc_k\|^2 > 0,\ \ \|\alpha_k^c v_k^c\| \neq 0,\ \ \text{and}\ \ \|c_k\| - \|c_k + \alpha_k^c J_k v_k^c\| > 0;
\eequation
otherwise, if $\|c_k\| = 0$, then $\|J_k^Tc_k\| = 0$ and it follows that $v_k = 0$ is the unique solution to~\eqref{prob.normal}.  Popular choices for computing a normal direction satisfying the aforementioned conditions include any of various Krylov subspace methods, such as the linear conjugate gradient (CG) method; see, e.g., \cite{NoceWrig06}.

Before describing the algorithm's procedure for computing the tangential direction, let us first introduce assumptions that the algorithm makes related to the stochastic gradients $\{g_k\}$ and symmetric matrices $\{H_k\}$ that it employs.

\bassumption\label{ass.g}
  There exists $M_g \in \R{}_{>0}$ such that, for all $k \in \N{}$, the stochastic gradient $g_k$ has the properties that $\E_k[g_k] = \nabla f_k$ and $\E_k[\|g_k - \nabla f_k\|^2] \leq M_g$, where $\E_k[\cdot]$ denotes expectation with respect to the distribution of $\omega$ $($recall~\eqref{prob.opt}$)$ conditioned on the event that $x_k$ is the primal iterate in iteration $k \in \N{}$.
\eassumption

Combining Assumption~\ref{ass.g} with Jensen's Inequality, it holds for all $k \in \N{}$ that
\bequation\label{jensens-M}
  \E_k[\|\nabla f_k - g_k\|] \leq \sqrt{\E_k[\|\nabla f_k - g_k\|^2]} \leq \sqrt{M_g}.
\eequation

\bassumption\label{ass.H}
  For all $k \in \N{}$, the matrix $H_k \in \Smbb^n$ is chosen independently from $g_k$.  In addition, there exist $M_H \in \R{}_{>0}$ and $\zeta \in\R{}_{>0}$ such that, for all $k \in \N{}$, it holds that $\|H_k\| \leq M_H$ and $u^TH_ku \geq \zeta\|u\|^2$ for all $u \in \Null(J_k)$.
\eassumption

For describing the tangential direction computation as it is performed in our algorithm, let us first describe what would be the computation of a tangential direction in a determinstic variant of our approach.  In particular, given $(x_k,y_k)$, $\nabla f_k$, a normal direction $v_k\in \Range(J_k^T)$, and $H_k$ satisfying Assumption~\ref{ass.H}, consider the subproblem
\bequation\label{prob.tangential_true}
  \min_{u\in\R{n}}\ (\nabla f_k + H_kv_k)^Tu + \tfrac{1}{2}u^TH_ku \ \  \text{s.t.}\ J_ku = 0,
\eequation
which has the unique solution $\uTrue_k \in \Null(J_k)$ that satisfies, for some $\deltaTrue_k \in \R{m}$,
\bequation\label{eq.linear_system_true}
  \bbmatrix H_k & J_k^T \\ J_k & 0 \ebmatrix \bbmatrix \uTrue_k \\ \deltaTrue_k \ebmatrix = -\bbmatrix \nabla f_k + H_kv_k + J_k^Ty_k \\ 0 \ebmatrix.
\eequation
This allows us to define, for purposes of our analysis only, the \emph{true and exact} primal-dual search direction (conditioned on $x_k$ being the $k$th iterate) as $(\dTrue_k,\deltaTrue_k)$, where $\dTrue_k := v_k + \uTrue_k$.  Since our algorithm only has access to a stochastic gradient estimate in each iteration, the corresponding \emph{exact} (but not \emph{true}) primal-dual search direction is given by $(d_{k,*},\delta_{k,*})$, where $d_{k,*} := v_k + u_{k,*}$ with $(u_{k,*},\delta_{k,*})$ satisfying
\bequation\label{eq.linear_system_exact}
  \bbmatrix H_k & J_k^T \\ J_k & 0 \ebmatrix \bbmatrix u_{k,*} \\ \delta_{k,*} \ebmatrix = -\bbmatrix g_k + H_kv_k + J_k^Ty_k \\ 0 \ebmatrix.
\eequation

Our algorithm, to avoid having to form or factor the matrix in \eqref{eq.linear_system_exact} in order to solve the system exactly, computes a tangential direction by computing $(u_k,\delta_k)$ through iterative linear algebra techniques applied to the symmetric indefinite system~\eqref{eq.linear_system_exact}.  In particular, our algorithm computes $(u_k,\delta_k)$ such that the full primal search direction $d_k := v_k + u_k$, dual search direction $\delta_k$, and residual defined by
\bequation\label{eq.linear_system}
  \bbmatrix \rho_k \\ r_k \ebmatrix := \bbmatrix H_k & J_k^T \\ J_k & 0 \ebmatrix \bbmatrix u_k \\ \delta_k \ebmatrix + \bbmatrix g_k + H_kv_k + J_k^Ty_k \\ 0 \ebmatrix
\eequation
satisfy at least one of a couple sets of conditions.  In the remainder of this subsection, we describe the sets of conditions that the algorithm employs to determine what constitutes an acceptable search direction (and corresponding pair of residuals).

In the deterministic setting, line search SQO methods commonly combine the search direction with an updating strategy for the merit parameter in a manner that ensures that the computed search direction is one of sufficient descent for the merit function.  The required descent condition is guaranteed to be satisfied by choosing the merit parameter to be sufficiently small so that the reduction in a model of the merit function (recall \eqref{eq.model_reduction}) is sufficiently large; see, e.g., \cite[Lemma~3.1]{ByrdCurtNoce08}.  Following such an approach, our algorithm requires that $(u_k,\delta_k)$ (yielding $d_k := v_k + u_k$) be computed and the merit parameter~$\tau$ be set such that the model reduction condition
\bequation\label{eq.model_reduction_condition}
  \Delta l(x_k,\tau,g_k,v_k+u_k) \geq \sigma_u\tau\max\{u_k^TH_ku,\epsilon_u\|u_k\|^2\} + \sigma_c (\|c_k\| - \|c_k + J_kv_k\|)
\eequation
holds for some user-defined $\sigma_u\in (0,1)$, $\epsilon_u \in (0,\zeta)$, and $\sigma_c\in (0,1)$.  (The particular value for the merit parameter $\tau$ for which the inequality \eqref{eq.model_reduction_condition} is required to hold depends on one of two different situations, as described below.)

Condition \eqref{eq.model_reduction_condition} plays a central role in the conditions that we require $(u_k,\delta_k)$ to satisfy.  We define these in the context of \emph{termination tests}, since they dictate conditions that, once satisfied, can cause termination of an iterative linear system solver applied to \eqref{eq.linear_system_exact}.  (The tests are inspired by the \emph{sufficient merit approximation reduction termination tests} developed in \cite{ByrdCurtNoce08,ByrdCurtNoce10,CurtNoceWach09} for the deterministic setting.)  Our first termination test states that an inexact solution of this linear system is acceptable if the model reduction condition~\eqref{eq.model_reduction_condition} is satisfied with the current merit parameter value (i.e., $\tau \equiv \tau_k \gets \tau_{k-1}$), the norms of the residual vectors satisfy certain upper bounds, and either the tangential direction is sufficiently small in norm compared to the normal direction or the tangential direction is one of sufficiently positive curvature for $H_k$ and yields a sufficiently small objective value for \eqref{prob.tangential_true} (with~$g_k$ in place of $\nabla f_k$).  The test makes use of a sequence $\{\beta_k\}$ that will also play a critical role in our step size selection scheme that is described in the next subsection.
\smallskip

\begin{tcolorbox}[colback=white]
\textbf{Termination Test 1.} 
  Given $\kappa\in (0,1)$, $\beta_k \in (0,1]$, $\kappa_\rho \in \R{}_{>0}$, $\kappa_r \in \R{}_{>0}$, $\kappa_u \in \R{}_{> 0}$, $\epsilon_u \in (0,\zeta)$, $\kappa_v \in \R{}_{> 0}$, $\sigma_u \in (0,1)$, $\sigma_c \in (0,1)$, and $v_k\in\Range(J_k^T)$ computed to satisfy~\eqref{eq.Cauchy_decrease}, the pair $(u_k,\delta_k)$ satisfies Termination Test~1 if, with the pair $(\rho_k,r_k)$ defined in~\eqref{eq.linear_system}, it holds that
  \bequation\label{eq.dual_residual_condition}
    \|\rho_k\| \leq \kappa \min \left\{ \left\|\bbmatrix g_k + J_k^T(y_k+\delta_k) \\ c_k\ebmatrix \right\|, \left\|\bbmatrix g_{k-1} + J_{k-1}^Ty_k \\ c_{k-1}\ebmatrix \right\| \right\};
  \eequation
  \bequation\label{eq.pd_residual_condition}
    \|\rho_k\| \leq \kappa_\rho\beta_k\ \ \text{and} \ \ \|r_k\| \leq \kappa_r \beta_k;
  \eequation
  \bequation\label{eq.tangential_component_condition_1}
    \|u_k\| \leq \kappa_u\|v_k\|\ \ \text{or} \ \ \left\{ \baligned u_k^TH_ku_k &\geq \epsilon_u\|u_k\|^2 \ \ \text{and} \\ (g_k + H_kv_k)^Tu_k + \tfrac{1}{2}u_k^TH_ku_k &\leq \kappa_v\|v_k\| \ealigned \right\};
  \eequation
  and \eqref{eq.model_reduction_condition} is satisfied with $\tau \equiv \tau_{k-1}$.  (In this case, the algorithm will set $\tau_k \gets \tau_{k-1}$ so that \eqref{eq.model_reduction_condition} holds with $\tau \equiv \tau_k$.)
\end{tcolorbox}

Termination Test~1 cannot be enforced in every iteration in a run of the algorithm, even in the deterministic setting, since there may exist points in the search space at which all of the conditions required in the test cannot be satisfied simultaneously, even if the linear system \eqref{eq.linear_system_exact} is solved to arbitrary accuracy.  In short, the algorithm needs to allow for the computation of a search direction for which \eqref{eq.model_reduction_condition} can only be satisfied with a decrease of the merit parameter.  That said, the algorithm needs to be careful in terms of the situations in which such a decrease is allowed to occur, or else the merit parameter sequence may behave in a manner that ruins a convergence guarantee for solving the original constrained optimization problem.  For our algorithm, we employ the following termination test for this situation.
\smallskip

\begin{tcolorbox}[colback=white]
\textbf{Termination Test 2.} 
  Given $\kappa \in (0,1)$, $\beta_k \in (0,1]$, $\kappa_\rho \in \R{}_{>0}$, $\kappa_r \in \R{}_{>0}$, $\kappa_u \in \R{}_{>0}$, $\epsilon_u \in (0,\zeta)$, $\kappa_v \in \R{}_{>0}$, $\epsilon_r \in (\sigma_c,1)$, and $v_k\in\Range(J_k^T)$ computed to satisfy~\eqref{eq.Cauchy_decrease}, the pair $(u_k,\delta_k)$ satisfies Termination Test~2 if, with the pair $(\rho_k,r_k)$ defined in~\eqref{eq.linear_system}, the conditions \eqref{eq.dual_residual_condition}--\eqref{eq.tangential_component_condition_1} hold along with
  \begin{align}\label{eq:TT2}
    \|c_k\| - \|c_k + J_kv_k + r_k\| \geq \epsilon_r(\|c_k\| - \|c_k + J_kv_k\|) > 0.
  \end{align}
  (In this case, for user-defined $\epsilon_{\tau}\in (0,1)$, the algorithm will set
  \bequation\label{eq:tau_update}
    \tau_k \gets \bcases \tau_{k-1} \quad &\text{if }\tau_{k-1}\leq \tautrial_k \\ \min\{(1-\epsilon_{\tau})\tau_{k-1},\tautrial_k\}\quad &\text{otherwise,} \ecases
  \eequation
  where
  \bequation\label{tautrail-def}
    \tautrial_k \gets \bcases \infty  &\text{if }g_k^Td_k + \max\{u_k^TH_ku_k,\epsilon_u\|u_k\|^2\} \leq 0 \\ \tfrac{(1-\tfrac{\sigma_c}{\epsilon_r})(\|c_k\| - \|c_k + J_kv_k + r_k\|)}{g_k^Td_k + \max\{u_k^TH_ku_k,\epsilon_u\|u_k\|^2\}}  &\text{otherwise,}  \ecases
  \eequation
  so~\eqref{eq.model_reduction_condition} is satisfied with $\tau \equiv \tau_k$.  See Lemma~\ref{lm:model_reduction_true} for a proof.)
\end{tcolorbox}

In Lemma~\ref{lem.direction_well_posed}, we show under a loose assumption about the behavior of the iterative linear system solver and a practical assumption about the algorithm iterates that, for all $k \in \N{}$, the algorithm can compute a pair $(u_k,\delta_k)$ satisfying at least one of Termination Test 1 or 2.  Therefore, the index of each iteration of our method is contained in one of two index sets, namely,
\bequationNN
  \baligned
    \Tone &:= \{k\in\N{}: (u_k,\delta_k) \ \text{satisfies Termination Test 1} \} \ \ \text{or} \\
    \Ttwo &:= \{k\in\N{}: (u_k,\delta_k) \ \text{satisfies Termination Test 2, but not Termination Test 1} \}.
  \ealigned
\eequationNN

\subsection{Computing a step size}

Upon the computation of $d_k \gets v_k + u_k$, our algorithm computes a positive step size $\alpha_k$ to determine $x_{k+1}$.  Given positive Lipschitz constants $L$ and~$\Gamma$ (recall Assumption~\ref{ass.main}), it follows for all $\alpha \in \R{}_{>0}$ that
\bequation\label{lip-consts}
  \baligned
    f(x_k + \alpha d_k) &\leq f_k + \alpha \nabla f_k^Td_k + \tfrac{1}{2}L\alpha^2\|d_k\|^2 \\
\text{and}\ \ \|c(x_k + \alpha d_k)\| &\leq \|c_k + \alpha J_kd_k\| + \tfrac{1}{2}\Gamma\alpha^2\|d_k\|^2.
  \ealigned
\eequation
Combining these inequalities with the definitions \eqref{eq.merit} and \eqref{eq.model_reduction}, the triangle inequality, and the definition of $r_k$ in \eqref{eq.linear_system}, one finds that
\bequation\label{eq.merit_reduction_bound}
  \baligned
    &\ \phi(x_k+\alpha d_k,\tau_k) - \phi(x_k,\tau_k) \\
   =&\ \tau_k f(x_k+\alpha d_k) - \tau_k f_k + \|c(x_k + \alpha d_k)\| - \|c_k\| \\
\leq&\ \alpha \tau_k \nabla f_k^Td_k + \|c_k + \alpha J_kd_k\| - \|c_k\| + \tfrac{1}{2}(\tau_k L + \Gamma)\alpha^2\|d_k\|^2 \\
\leq&\ \alpha \tau_k \nabla f_k^Td_k + (|1-\alpha| - 1)\|c_k\| + \alpha\|c_k + J_kd_k\| + \tfrac{1}{2}(\tau_k L + \Gamma)\alpha^2\|d_k\|^2 \\
   =&\ -\alpha\Delta l(x_k,\tau_k,\nabla f_k,d_k) + (|1-\alpha| - (1-\alpha))\|c_k\| + \tfrac{1}{2}(\tau_k L + \Gamma)\alpha^2\|d_k\|^2.
  \ealigned
\eequation
This derivation provides a convex piecewise-quadratic upper-bounding function for the change in the merit function corresponding to a step from $x_k$ to $x_k + \alpha d_k$.  Given user-defined $\eta\in(0,1)$ and the aforementioned sequence $\{\beta_k\}\subset(0,1]$, our algorithm's step size selection scheme makes use of the quantity 
\bequation\label{eq:stochastic_alpha_suff}
  \alphasuff_k := \min \left\{\tfrac{2(1-\eta)\beta_k\Delta l(x_k,\tau_k,g_k,d_k)}{(\tau_k L+\Gamma)\|d_k\|^2}, 1 \right\}.
\eequation   
The definition of $\alphasuff_k$ can be motivated as follows.  Its value, when $\beta_k = 1$, is the largest value on $[0,1]$ such that for all $\alpha\in[0,\alphasuff_k]$ the right-hand-side of~\eqref{eq.merit_reduction_bound} (with $\nabla f_k$ replaced by $g_k$) is less than or equal to $-\eta\alpha\Delta l(x_k,\tau_k,g_k,d_k)$.  Such an inequality is representative of one enforced in deterministic line search SQO methods.  Otherwise, with $\beta_k \in (0,1]$ introduced and not necessarily equal to 1, the value of $\alphasuff_k$ can be diminished over the course of the optimization process, which allows for step size control as is required for convergence guarantees for certain stochastic-gradient-based methods; see, e.g., \cite{BottCurtNoce18}.  The first term inside the min appearing in~\eqref{eq:stochastic_alpha_suff} is important for the convergence guarantees that we prove for our method, but it can behave erratically due to the algorithm's use of stochastic gradient estimates.  To account for this stochasticity, given user-defined $\epsilon_{\xi}\in (0,1)$, our algorithm defines
\bequation\label{eq:stochastic_ratio}
  \xitrial_k := \tfrac{\Delta l(x_k,\tau_k,g_k,d_k)}{\tau_k\|d_k\|^2}\ \ \text{and} \ \ \xi_k := 
\bcases \xi_{k-1}\qquad &\text{ if }\xi_{k-1}\leq \xitrial_k \!\!\\ 
\min\{(1-\epsilon_{\xi})\xi_{k-1},\xitrial_k\} &\text{ otherwise,} \ecases
\eequation
so that $\xi_k \leq \xitrial_k = \Delta l(x_k,\tau_k,g_k,d_k)/(\tau_k\|d_k\|^2)$ for all $k\in\N{}$. Combining this inequality with~\eqref{eq:stochastic_alpha_suff}, the monotonically nonincreasing behaviors of $\{\xi_k\}$ and $\{\tau_k\}$, and assuming that the sequence $\{\beta_k\}$ is chosen to satisfy
\bequation\label{betak-requirement}
  2(1-\eta)\beta_k\xi_{-1}\tau_{-1}/\Gamma \in (0,1]\ \ \text{for all} \ \ k\in\N{}
\eequation
where $\xi_{-1}$ and $\tau_{-1}$ initialize the sequences $\{\xi_k\}$ and $\{\tau_k\}$, one finds
\bequation\label{eq:stochastic_alpha_suff.new}
  \alphamin_k := \tfrac{2(1-\eta)\beta_k\xi_k\tau_k}{(\tau_k L+\Gamma)} \leq \min \left\{\tfrac{2(1-\eta)\beta_k\Delta l(x_k,\tau_k,g_k,d_k)}{(\tau_k L+\Gamma)\|d_k\|^2}, 1 \right\} \equiv \alphasuff_k.
\eequation 
The value $\alphamin_k$ serves as a minimum value (i.e., a lower bound) for our choice of step size.  In our analysis, we will also show that even though $\xi_k$ is stochastic for each $k\in\N{}$, the sequence $\{\xi_k\}$ is bounded away from zero deterministically (see Lemma~\ref{lem:xi-bound}).

Next, let us derive a maximum value (i.e., an upper bound) for our algorithm's choice of step size. Consider the strongly convex function $\varphi : \R{} \to \R{}$ defined by
\bequation\label{eq:stochastic_U_def}
  \baligned
    \varphi(\alpha) := &(\eta - 1)\alpha\beta_k\Delta l(x_k,\tau_k,g_k,d_k) + \|c_k + \alpha J_kd_k\| - \|c_k\| \\
&+ \alpha(\|c_k\| - \|c_k + J_kd_k\|) + \tfrac{1}{2}(\tau_k L+\Gamma)\alpha^2\|d_k\|^2.
  \ealigned
\eequation 
Notice that when $\beta_k = 1$, it holds that $\varphi(\alpha) \leq 0$ for all $\alpha \in \R{}_{\geq0}$ if and only if the quantity in the third row of~\eqref{eq.merit_reduction_bound} (with $\nabla f_k$ replaced by $g_k$) is less than or equal to $-\eta\alpha\Delta l(x_k,\tau_k,g_k,d_k)$.  Thus, following a similar argument as above, one can be motivated as to the fact that our algorithm never allows a step size larger than
\bequation\label{eq.alpha_phi}
  \alphavarphi_k := \max\{ \alpha \in \R{}_{\geq0} : \varphi(\alpha) \leq 0\}.
\eequation
Finally, again to mitigate adverse affects caused by the use of stochastic gradient estimates, our algorithm employs the maximum step size
\bequation\label{def-alphamax}
  \alphamax_k := \min \{ \alphavarphi_k, \alphamin_k + \theta\beta_k^2\},
\eequation
where $\theta\in\R{}_{\geq 0}$ is user-defined.  Overall, our algorithm allows any step size with $\alpha_k \in [\alphamin_k,\alphamax_k]$.  Lemma~\ref{lm:stepsize_well_defined} shows that this interval is nonempty.

\subsection{Updating the primal-dual iterate}

In the primal space, our algorithm employs the iterate update $x_{k+1} \gets x_k + \alpha_k d_k$.  However, in the dual space, it allows additional flexibility; in particular, the algorithm allows any $y_{k+1}$ such that
\bequation\label{eq:stochastic_dual_update}
  \|g_k + J_k^Ty_{k+1}\|\leq \|g_k + J_k^T(y_k + \delta_k)\|.
\eequation
Clearly, choosing $y_{k+1} \gets y_k + \delta_k$ is one particular option satisfying \eqref{eq:stochastic_dual_update}, although other choices such as least-squares multipliers could also be used.

\balgorithm[ht]
  \caption{Stochastic Inexact Sequential Quadratic Optimization (SISQO)}
  \label{alg:stochastic_sqp_adaptive}
  \balgorithmic[1]
    \Require initial values $(x_0,y_0,\tau_{-1},\xi_{-1}) \in \R{n} \times \R{m} \times \R{}_{>0} \times \R{}_{>0}$; Lipschitz constants $(L,\Gamma) \in \R{}_{>0} \times \R{}_{>0}$ satisfying Assumption~\ref{ass.main}; $\epsilon_c\in (0,1]$; $\epsilon_u\in (0,\zeta)$; $\{\sigma_u, \sigma_c, \kappa, \epsilon_{\tau}, \epsilon_{\xi}, \eta\} \subset (0,1)$; $\{\kappa_\rho, \kappa_r, \kappa_u, \kappa_v, \theta\} \subset \R{}_{>0}$; $\epsilon_r \in (\sigma_c,1)$
    \For{\textbf{all} $k \in \N{}$}
      \State choose $\beta_k \in (0,1]$ satisfying~\eqref{betak-requirement} and $H_k$ satisfying Assumption~\ref{ass.H}
      \State compute $v_k\in\Range(J_k^T)$ satisfying \eqref{eq.Cauchy_decrease}
      \State generate $g_k$ satisfying Assumption~\ref{ass.g}
      \State compute $(u_k,\delta_k)$ satisfying at least one of Termination Tests 1 or 2 \label{line:stochastic_iterative_solver}
      \If{Termination Test 1 is satisfied} 
        \State set $\tautrial_k \gets \infty$ and $\tau_k \gets \tau_{k-1}$ \Comment{$k\in\Tone$}
      \Else{ (Termination Test 2 is satisfied)}
        \State set $\tautrial_k$ and $\tau_k$ by \eqref{eq:tau_update}--\eqref{tautrail-def} \label{line:stochastic_termination_test_2} \Comment{$k\in\Ttwo$}
      \EndIf   
      \State set $d_k \gets v_k + u_k$
      \State compute $\xi_k$ and $\xitrial_k$ by \eqref{eq:stochastic_ratio} \label{line:stochastic_update_xi}
      \State choose $\alpha_k\in [\alphamin_k,\alphamax_k]$ using the definitions in~\eqref{eq:stochastic_alpha_suff.new} and \eqref{def-alphamax} \label{line:stochastic_choose_stepsize}
      \State set $x_{k+1} \gets x_k + \alpha_k d_k$ and choose $y_{k+1}$ satisfying \eqref{eq:stochastic_dual_update}
	\EndFor
  \ealgorithmic
\ealgorithm

\section{Analysis}\label{sec.analysis}

Our analysis is presented in three parts.  In Section~\ref{sec.wp}, we show that Algorithm~\ref{alg:stochastic_sqp_adaptive} is well posed.  Then, in Section~\ref{sec.general}, we prove general lemmas about the behavior of our algorithm.  Finally, in Section~\ref{sec.convergence}, we prove convergence properties in expectation for the iterate sequence generated by the algorithm under an assumption about the behavior of the merit parameter sequence.  The assumption employed can be justified using the same arguments as in \cite{BeraCurtRobiZhou21}, as explained in Section~\ref{sec.convergence}.

\subsection{Well-posedness}\label{sec.wp}

Our aim in this subsection is to prove that during each iteration of Algorithm~\ref{alg:stochastic_sqp_adaptive}, each step of the algorithm can be performed in a manner that terminates finitely.  Along the way, we also establish useful properties of quantities computed by the algorithm.  We make the following reasonable assumption concerning the behavior of the iterative linear system solver employed by the algorithm for the tangential direction and dual step computation.

\bassumption\label{ass:sp-solver}
  For all $k\in\N{}$, the iterative linear system solver employed in line~\ref{line:stochastic_iterative_solver} generates a sequence $\{(u_{k,t},\delta_{k,t},\rho_{k,t},r_{k,t})\}_{t\in\N{}}$ satisfying
  \bequation\label{eq:stochastic_system_iterative}
    \bbmatrix \rho_{k,t} \\ r_{k,t} \ebmatrix
= \bbmatrix H_k & J_k^T \\ J_k & 0 \ebmatrix \bbmatrix u_{k,t} \\ \delta_{k,t} \ebmatrix +\bbmatrix g_k + H_kv_k + J_k^T y_k \\ 0 \ebmatrix\ \ \text{for all}\ \ t \in \N{}
  \eequation
  such that $\lim_{t\to\infty} \| (u_{k,t},\delta_{k,t},\rho_{k,t},r_{k,t})- (u_{k,*},\delta_{k,*},0,0) \| = 0$, where $(u_{k,*},\delta_{k,*})$ is the unique solution to the linear system defined in~\eqref{eq.linear_system_exact}.
\eassumption

We also make the following assumption concerning the algorithm iterates and corresponding stochastic gradient estimates computed in each iteration.

\bassumption\label{ass:nonzero}
  For all $k\in\N{}$, it holds that $c_k \neq 0$ or $g_k \notin \Range(J_k^T)$.
\eassumption

We justify Assumption~\ref{ass:nonzero} in the following manner.  In the deterministic setting, the algorithm encounters a point $x_k$ such that $c_k = 0$ and $\nabla f_k \in \Range(J_k^T)$ if and only if there exists $y_k$ such that \eqref{eq.stationary} holds for $(x,y) \equiv (x_k,y_k)$, i.e., the point $(x_k,y_k)$ is first-order stationary for problem \eqref{prob.opt}.  In such a scenario, it is reasonable to require that an exact solution of \eqref{eq.linear_system_exact} is computed, or at least a sufficiently accurate solution of the system is computed such that a practical termination condition for \eqref{eq.linear_system_exact} is triggered and the algorithm terminates.  In the stochastic setting, the algorithm encounters $c_k = 0$ and $g_k \in \Range(J_k^T)$ if and only if $x_k$ is \emph{exactly} feasible and the stochastic gradient lies \emph{exactly} in the range space of~$J_k^T$.  Since $g_k$ is a \emph{stochastic} gradient, we contend that it is unlikely that it will lie \emph{exactly} in $\Range(J_k^T)$ except in special circumstances.  Thus, for simplicity in our analysis, we impose Assumption~\ref{ass:nonzero} throughout this section.  (If Assumption~\ref{ass:nonzero} were not to hold, then one of the following could be employed in a practical implementation: (i) if a sufficiently accurate solution of \eqref{eq.linear_system_exact} does not satisfy either Termination Test 1 or 2, then a new stochastic gradient could be sampled, perhaps following a procedure to ensure that if multiple new stochastic gradients are computed, then each is computed with lower variance, or (ii) random (e.g., Gaussian) noise could be added to $g_k$ for all $k \in \N{}$ so that Assumption~\ref{ass:nonzero} holds with probability one in all iterations, in which case the convergence result that we prove will hold with probability one.)

We can now show that the search direction computation is well posed.

\blemma\label{lem.direction_well_posed}
  For all $k\in\N{}$, the iterative linear system solver computes $(u_k,\delta_k)$ satisfying at least one of Termination Test~1 or 2 in a finite number of iterations.
\elemma
\bproof
  We prove the result by considering two cases.
  
  \textbf{Case 1:} $\|c_k\| > 0$. For this case, we show that $(u_k,\delta_k) \equiv (u_{k,t},\delta_{k,t})$ satisfies Termination Test~2 for sufficiently large $t \in \N{}$.  Let us first observe that it follows from  Assumption~\ref{ass:sp-solver}, Assumption~\ref{ass:nonzero}, and the fact that $\beta_k \in (0,1]$ that both  \eqref{eq.dual_residual_condition} and~\eqref{eq.pd_residual_condition} hold with $(\rho_k,r_k) \equiv (\rho_{k,t},r_{k,t})$ for all sufficiently large $t \in \N{}$.

  Let us now show that \eqref{eq.tangential_component_condition_1} holds for all sufficiently large $t \in \N{}$.  Since $\|c_k\| > 0$, it follows under Assumption~\ref{ass.main} that $\|v_k\| > 0$.  If $\|u_{k,*}\| = 0$, then Assumption~\ref{ass:sp-solver} implies $\{\|u_{k,t}\|\}\to\|u_{k,*}\| = 0$, in which case it follows from $\kappa_u \in \R{}_{>0}$ that the former condition in~\eqref{eq.tangential_component_condition_1} holds with $u_k \equiv u_{k,t}$ for all sufficiently large $t \in \N{}$. On the other hand, if $\|u_{k,*}\| > 0$, then \eqref{eq.linear_system_exact} and Assumption~\ref{ass.H} imply
  \bequation\label{star-eq}
    \baligned
      u_{k,*}^T(g_k+H_k v_k) + \tfrac{1}{2} u_{k,*}^T H_k u_{k,*} &< u_{k,*}^T(g_k+H_k v_k) + u_{k,*}^T H_k u_{k,*} \\
      &= u_{k,*}^T(g_k+H_k v_k + H_k u_{k,*} + J_k^Ty_k) \\
      &= -u_{k,*}^T J_k^T\delta_{k,*} = - (J_k u_{k,*})^T \delta_{k,*} = 0.
    \ealigned
  \eequation
  Combining this inequality with the facts that $\epsilon_u \in (0,\zeta)$, $\kappa_v \in \R{}_{>0}$, and $\|v_k\| > 0$, it follows under Assumptions~\ref{ass.H} and~\ref{ass:sp-solver} that the latter set of conditions in~\eqref{eq.tangential_component_condition_1} holds with $u_k \equiv u_{k,t}$ for all sufficiently large $t \in \N{}$.

  Finally, let us show that \eqref{eq:TT2} holds for all sufficiently large $t \in \N{}$, which combined with the previous conclusions shows that Termination Test~2 is satisfied by $(u_k,\delta_k) \equiv (u_{k,t},\delta_{k,t})$ for all sufficiently large $t \in \N{}$.  By Assumption~\ref{ass:sp-solver}, \eqref{eq.Cauchy_decrease}, and the aforementioned fact that $\|v_k\| > 0$, it follows that
  \bequationNN
    \lim_{t\to\infty} (\|c_k\|-\|c_k + J_kv_k + r_{k,t}\|) = \|c_k\| - \|c_k + J_kv_k\| > 0,
  \eequationNN
  which shows that \eqref{eq:TT2} holds with $r_k \equiv r_{k,t}$ for all sufficiently large $t \in \N{}$, as desired.
  
  \textbf{Case 2:} $\|c_k\| = 0$. For this case, we show that $(u_k,\delta_k) \equiv (u_{k,t},\delta_{k,t})$ satisfies Termination Test~1 for all sufficiently large $t \in \N{}$.  First, recall that $\|c_k\| = 0$ implies that $\|v_k\| = 0$.  We also claim that $\|u_{k,*}\| > 0$.  To prove this by contradiction, suppose that $\|u_{k,*}\| = 0$.  Combining this with $v_k = 0$ and~\eqref{eq.linear_system_exact}, it follows that $g_k + J_k^T(y_k+\delta_{k,*}) = 0$, which with $c_k = 0$ violates Assumption~\ref{ass:nonzero}.  Thus, $\|u_{k,*}\| > 0$.

  Next, notice that the argument used in the beginning of Case 1 still applies in this case, which allows us to conclude that both \eqref{eq.dual_residual_condition} and~\eqref{eq.pd_residual_condition} hold with $(\rho_k,r_k) = (\rho_{k,t},r_{k,t})$ for all sufficiently large $t \in \N{}$.  Also,  since $\|u_{k,*}\| > 0$, the first inequality in~\eqref{star-eq} holds as a strict inequality, i.e., $u_{k,*}^T(g_k+H_k v_k) + \tfrac{1}{2} u_{k,*}^T H_k u_{k,*} < 0$. Combining this  inequality with  Assumption~\ref{ass:sp-solver}, Assumption~\ref{ass.H}, and $\epsilon_u\in(0,\zeta)$ allows us to deduce that the second set of conditions in~\eqref{eq.tangential_component_condition_1} holds with $u_k \equiv u_{k,t}$ for all sufficiently large $t \in \N{}$.  Next, from the fact that $\|v_k\|= 0$ and~\eqref{eq.linear_system_exact}, it follows that $J_kd_{k,*} = J_k(u_{k,*} + v_k) = 0$, which with Assumption~\ref{ass.H} and $\epsilon_u\in(0,\zeta)$ gives $u_{k,*}^T H_k u_{k,*} \geq \zeta \|u_{k,*}\|^2 > \epsilon_u\|u_{k,*}\|^2$, from which we deduce that $\max\{u_{k,*}^T H_k u_{k,*}, \epsilon_u\|u_{k,*}\|^2\} = u_{k,*}^T H_k u_{k,*} \geq \zeta \|u_{k,*}\|^2 > 0$.  Combining this inequality with $\|c_k\| = 0$, $\|v_k\| = 0$, $J_kd_{k,*} = J_kv_k = 0$, \eqref{eq.linear_system_exact}, and Assumption~\ref{ass.H} shows that
  \bequationNN
    \baligned
      &\ \Delta l(x_k,\tau_{k-1},g_k,d_{k,*}) = -\tau_{k-1}g_k^Td_{k,*} + \|c_k\| - \|c_k+J_k d_{k,*}\| = -\tau_{k-1} g_k^T u_{k,*} \\
     =&\ -\tau_{k-1}(-H_ku_{k,*} - H_kv_k - J_k^T(y_k + \delta_{k,*}))^Tu_{k,*} = \tau_{k-1}u_{k,*}^TH_ku_{k,*} \\
     >&\ \sigma_u\tau_{k-1}\max\{u_{k,*}^T H_k u_{k,*}, \epsilon_u\|u_{k,*}\|^2\} + \sigma_c(\|c_k\| - \|c_k + J_kv_k\|) > 0,
    \ealigned
  \eequationNN
  meaning that the sufficient decrease condition~\eqref{eq.model_reduction_condition} holds with $\tau \equiv \tau_{k-1}$ for all sufficiently large $t \in \N{}$.  In summary, we have shown that, for all sufficiently large $t \in \N{}$, the pair $(u_k,\delta_k) \equiv (u_{k,t},\delta_{k,t})$ will satisfy Termination Test~1, as desired.  
\eproof

Next, we prove that every search direction is nonzero in norm.

\blemma\label{lem:d-nonzero}
  For all $k\in\N{}$, it holds that $\|d_k\| > 0$.
\elemma
\bproof
  For a proof by contradiction, suppose that $\|d_k\| = 0$.  From this fact, $d_k = v_k+u_k$, and~\eqref{eq.linear_system}, it follows that $\rho_k = g_k+J_k^T(y_k+\delta_k) + H_k(v_k+u_k) = g_k + J_k^T(y_k+\delta_k)$. If $\|c_k\| = 0$, then this value for $\rho_k$ shows that the inequality in~\eqref{eq.dual_residual_condition} cannot hold, meaning that $(u_k,\delta_k)$ cannot satisfy Termination Test 1 or 2, which contradicts Lemma~\ref{lem.direction_well_posed}.  Hence, the only possibility is that $\|c_k\| > 0$, which we shall assume for the remainder of the proof.

  Notice from $\|d_k\| = 0$, $d_k = v_k+u_k$, and $r_k = J_ku_k$, it follows that $\|c_k\| - \|c_k+J_kv_k + r_k\| = \|c_k\| - \|c_k+J_kd_k\| = 0$, meaning that~\eqref{eq:TT2} is not satisfied; thus, $(u_k,\delta_k)$ does not satisfy Termination Test 2.  Also, observe from $\|v_k\| > 0$ (which follows from $\|c_k\| > 0$ and Assumption~\ref{ass.main}), $\|d_k\| = 0$, and \eqref{eq.Cauchy_decrease} that $\Delta l(x_k,\tau_k,g_k,d_k) = 0 < \sigma_u\tau_{k-1}\max\{u_k^TH_ku_k,\epsilon_u\|u_k\|^2\} + \sigma_c (\|c_k\| - \|c_k + J_kv_k\|)$, meaning that~\eqref{eq.model_reduction_condition} is not satisfied with $\tau = \tau_{k-1}$; thus, $(u_k,\delta_k)$ does not satisfy Termination Test 1.  Overall, we have reached a contradiction to Lemma~\ref{lem.direction_well_posed}, and since we have reached a contradiction in all cases, the original supposition that $\|d_k\| = 0$ cannot be true.
\eproof

We now show that our update strategy for the merit parameter sequence ensures that the model reduction condition \eqref{eq.model_reduction_condition} always holds for $\tau \equiv \tau_k$.  We also show another important property of the sequence $\{\tau_k\}$.

\blemma\label{lm:model_reduction_true}
  For all $k\in\N{}$, the inequality in  \eqref{eq.model_reduction_condition} holds with $\tau \equiv \tau_k$.  In addition, for all $k \in \N{}$ such that $\tau_{k+1} < \tau_k$, it holds that $\tau_{k+1} \leq (1-\epsilon_{\tau})\tau_k$.
\elemma
\bproof
  The desired conclusion follows for $k\in \Tone$ due to the manner in which Termination Test~1 is defined and the fact that the algorithm sets $\tau_k \gets \tau_{k-1}$ for all $k \in \Tone$.  Hence, let us proceed under the assumption that $k\in\Ttwo$.  The inequality in  \eqref{eq.model_reduction_condition} holds for $\tau \equiv \tau_k$ with $d_k = v_k+u_k$ if and only if
  \bequationNN
    \tau_k (g_k^Td_k + \sigma_u \max \{u_k^TH_ku_k,\epsilon_u\|u_k\|^2 \} ) \leq \|c_k\| - \|c_k + J_kd_k\| - \sigma_c(\|c_k\| - \|c_k + J_kv_k\|).
  \eequationNN
  We now proceed to show that this inequality holds by considering two cases.
  
  \textbf{Case 1:} $g_k^Td_k + \max\{u_k^TH_ku_k,\epsilon_u\|u_k\|^2\} \leq 0$.  In this case, the algorithm sets $\tau_k \gets \tau_{k-1}$.  Combining this with~\eqref{eq:TT2}, $J_ku_k = r_k$, and $\epsilon_r\in(\sigma_c,1)$ yields
  \bequationNN
    \baligned
      \tau_k (g_k^Td_k &+ \sigma_u \max\{u_k^TH_ku_k,\epsilon_u\|u_k\|^2 \} ) \leq \tau_k (g_k^Td_k + \max \{u_k^TH_ku_k,\epsilon_u\|u_k\|^2 \} ) \\
\leq 0
      &\leq \|c_k\| - \|c_k + J_kd_k\| - \epsilon_r(\|c_k\| - \|c_k + J_kv_k\|) \\
      &< \|c_k\| - \|c_k + J_kd_k\| - \sigma_c(\|c_k\| - \|c_k + J_kv_k\|),
    \ealigned
  \eequationNN
  which establishes the desired inequality.

  \textbf{Case 2:} $g_k^Td_k + \max\{u_k^TH_ku_k,\epsilon_u\|u_k\|^2\} > 0$. The update~\eqref{eq:tau_update} yields $\tau_k \leq \tautrial_k$, which combined with \eqref{eq:TT2}, \eqref{tautrail-def}, $J_ku_k = r_k$, and $\epsilon_r\in(\sigma_c,1)$ yields
  \bequationNN
    \baligned
      &\ \tau_k (g_k^Td_k + \sigma_u \max \{u_k^TH_ku_k,\epsilon_u\|u_k\|^2 \} ) \leq \tau_k (g_k^Td_k + \max \{u_k^TH_ku_k,\epsilon_u\|u_k\|^2 \} ) \\
  \leq&\ (1-\tfrac{\sigma_c}{\epsilon_r})(\|c_k\| - \|c_k + J_kd_k\|) \leq \|c_k\| - \|c_k + J_kd_k\| - \sigma_c(\|c_k\| - \|c_k + J_kv_k\|),
    \ealigned
  \eequationNN
  as desired. Moreover, from \eqref{eq:tau_update}, we have $\tau_{k+1}\leq (1-\epsilon_{\tau})\tau_k$ whenever $\tau_{k+1} < \tau_k$.
\eproof

We conclude this subsection by showing that the interval defining our step size selection scheme, i.e., $[\alphamin_k,\alphamax_k]$, is positive and nonempty for all $k \in \N{}$.  We also show a useful property of the computed step size that is needed in our analysis.

\blemma\label{lm:stepsize_well_defined}
  For all $k\in\N{}$, it holds that $0 < \alphamin_k \leq \alphasuff_k \leq \alphavarphi_k$ and $0 < \alphamin_k \leq \alphamax_k$.  In addition, for all $k \in \N{}$, it holds that $\varphi(\alpha_k) \leq 0$.
\elemma
\bproof
  It follows from \eqref{eq:stochastic_alpha_suff.new} and the fact that $\{\beta_k\}$, $\{\xi_k\}$, and $\{\tau_k\}$ are positive sequences that $\alphamin_k > 0$ for all $k \in \N{}$.  Hence, considering \eqref{eq:stochastic_alpha_suff.new} and \eqref{def-alphamax}, to prove that $0 < \alphamin_k \leq \alphasuff_k \leq \alphavarphi_k$ and $0 < \alphamin_k \leq \alphamax_k$ for all $k \in \N{}$, it is sufficient to show that $\alphasuff_k \leq \alphavarphi_k$ for all $k \in \N{}$.  Consider arbitrary $k \in \N{}$.  Since $\alphavarphi_k \geq 0$ by construction and $\alphasuff_k \geq 0$ as a consequence of Lemmas~\ref{lem:d-nonzero} and \ref{lm:model_reduction_true}, the inequality holds trivially if $\alphasuff_k = 0$.  Hence, we may proceed under the assumption that $\alphasuff_k > 0$.  Moreover, one finds from the definition of $\alphavarphi_k$ in \eqref{eq.alpha_phi} that to establish $\alphasuff_k \leq \alphavarphi_k$ it is sufficient to show that $\varphi(\alphasuff_k) \leq 0$. We consider two cases based on the min in~\eqref{eq:stochastic_alpha_suff}.  First, suppose that $\alphasuff_k = 1 \leq \tfrac{2(1-\eta)\beta_k\Delta l(x_k,\tau_k,g_k,d_k)}{(\tau_k L+\Gamma)\|d_k\|^2}$, which with \eqref{eq:stochastic_U_def} shows that
  \bequationNN
    \baligned
      \varphi(\alphasuff_k) 
        &= (\eta-1)\beta_k\Delta l(x_k,\tau_k,g_k,d_k) + \tfrac{1}{2}(\tau_kL+\Gamma)\|d_k\|^2 \\
        &\leq (\eta-1)\beta_k\Delta l(x_k,\tau_k,g_k,d_k) + (1-\eta)\beta_k\Delta l(x_k,\tau_k,g_k,d_k) = 0,
    \ealigned
  \eequationNN 
  as desired. Second, suppose $\alphasuff_k = \tfrac{2(1-\eta)\beta_k\Delta l(x_k,\tau_k,g_k,d_k)}{(\tau_k L+\Gamma)\|d_k\|^2} < 1$.  For this case, it follows from  \eqref{eq:stochastic_U_def}, $\alphasuff_k \in (0,1]$, and the triangle inequality that
  \bequationNN
    \baligned
      \varphi(\alphasuff_k) 
        &= (\eta-1)\alphasuff_k\beta_k\Delta l(x_k,\tau_k,g_k,d_k) + (1-\eta)\alphasuff_k\beta_k\Delta l(x_k,\tau_k,g_k,d_k) \\
        &\phantom{=i} +\|c_k + \alphasuff_k J_kd_k\| - \alphasuff_k\|c_k + J_kd_k\| + (\alphasuff_k-1)\|c_k\| \\
        &\leq \|(1-\alphasuff_k)c_k\| + (\alphasuff_k-1)\|c_k\| = 0.
    \ealigned
  \eequationNN
  Overall, $\alphasuff_k \leq \alphavarphi_k$ since, in both cases above, we proved that $\varphi(\alphasuff_k) \leq 0$.
  
  Finally, let us show that $\varphi(\alpha_k) \leq 0$ for all $k \in \N{}$.  By \eqref{eq.model_reduction} and \eqref{eq:stochastic_U_def}, one finds (as previously mentioned) that $\varphi$ is strongly convex.  In addition, one finds that $\varphi(0) = \varphi(\alphavarphi_k) = 0$, where $\alphavarphi_k \in \R{}_{>0}$ due to the first part of this lemma.  Along with the fact that $0 < \alphamin_k \leq \alpha_k \leq \alphamax_k \leq \alphavarphi_k$, it follows that $\varphi(\alpha_k) \leq 0$, as desired.
\eproof

\subsection{General results}\label{sec.general}

Our aim in this subsection is to prove general results about the behavior of quantities generated by Algorithm~\ref{alg:stochastic_sqp_adaptive}.  For our purposes here, we make the following assumption about the dual and residual sequences.

\bassumption\label{ass:bd-residual}
  The dual iterate sequence $\{y_k\}$ and residual sequence $\{(\rho_k,r_k)\}$ $($recall~\eqref{eq.linear_system}$)$ generated by Algorithm~\ref{alg:stochastic_sqp_adaptive} are bounded in norm.
\eassumption

We note that under Assumption~\ref{ass.main} and Assumption~\ref{ass.H}, this additional assumption is mild; it should hold as long as any reasonable iterative solver is applied to \eqref{eq.linear_system_exact} in each iteration of a run of the algorithm.

The next lemma gives a lower bound on 
$\|c_k\| - \|c_k + J_kv_k\|$ relative to $\|c_k\|$.

\blemma\label{lm:feasibility_decrease}
  There exists $\omega_1\in\R{}_{>0}$ such that, for all $k\in\N{}$, it holds that
  \bequationNN
    \|c_k\| - \|c_k + J_kv_k\| \geq \omega_1 \|c_k\|.
  \eequationNN
\elemma
\bproof
  This result follows as in \cite[Lemma~3.5]{CurtNoceWach09}, but with small straightforward modifications to account for the fact that, in our analysis here, the singular values of $\{J_k\}$ are bounded away from zero as a consequence of Assumption~\ref{ass.main}.
\eproof

The next lemma shows that $\|v_k\|$ is of the same order as $\|c_k\|$.

\blemma\label{lm:vk_bounded}
  There exists $\{\omega_2,\omega_3\}\subset\R{}_{>0}$ such that, for all $k\in\N{}$, it holds that
  \bequation\label{eq:vk-bounded}
    \omega_2\|c_k\|\leq \|v_k\|\leq \omega_3\|c_k\|.
  \eequation
\elemma
\bproof
  Observe that Assumption~\ref{ass.main} ensures the existence of $\lambdamin\in\R{}_{>0}$ such that $J_kJ_k^T \succeq \lambdamin I$ for all $k\in\N{}$.  We now prove each desired inequality.  First, consider the former inequality in~\eqref{eq:vk-bounded}.  Since this inequality  holds trivially whenever $\|c_k\| = 0$, let us proceed under the assumption that $\|c_k\| > 0$.  One finds
  \begin{align*}
    \|c_k\|^2 - \|c_k+\alpha_k^c J_k v_k^c\|^2 
      &= (\|c_k\| - \|c_k+\alpha_k^c J_k v_k^c\|)(\|c_k\| + \|c_k+\alpha_k^c J_k v_k^c\|) \\ &\leq 2\|c_k\|(\|c_k\| - \|c_k+\alpha_k^c J_k v_k^c\|).
  \end{align*}
  It follows from this inequality, the triangle inequality, and \eqref{eq.Cauchy_decrease} that
  \bequationNN
    \baligned
      &\ \ \ \ \|J_k\|\|v_k\| \geq \|J_kv_k\| \geq \|c_k\| - \|c_k + J_kv_k\| \geq \epsilon_c(\|c_k\| - \|c_k + \alpha_k^c J_k v_k^c\|) \\
      &\geq \tfrac{\epsilon_c}{2\|c_k\|}(\|c_k\|^2 - \|c_k + \alpha_k^c J_k v_k^c\|^2) = \tfrac{\epsilon_c}{2\|c_k\|}(-2\alpha_k^c c_k^T J_k v_k^c - (\alpha_k^c)^2\|J_kv_k^c\|^2).
    \ealigned
  \eequationNN
  Substituting in for the value of $\alpha_k^c$ (recall \eqref{eq.alphac}), then substituting $v_k^c = -J_k^T c_k$ and simplifying shows that $\|J_k\|\|v_k\| \geq ( \tfrac{\epsilon_c}{2\|c_k\|}) \alpha_k^c\|J_k^Tc_k\|^2$.  Again substituting the value of $\alpha_k^c$ and using the definition of $\lambdamin$, it follows that
  \bequationNN
    \|J_k\|\|v_k\| \geq \tfrac{\epsilon_c \|J_k^Tc_k\|^4}{2\|c_k\|\|J_kJ_k^Tc_k\|^2}
\geq \tfrac{\epsilon_c\lambdamin^2\|c_k\|^4}{2\|c_k\|\|J_k^TJ_k\|^2\|c_k\|^2} 
= \tfrac{\epsilon_c\lambdamin^2}{2\|J_k^TJ_k\|^2}\|c_k\|.  
  \eequationNN
  It follows from this inequality and Assumption~\ref{ass.main} that there exists $\omega_2\in\R{}_{>0}$ such that the former inequality in~\eqref{eq:vk-bounded} holds, as desired.
  
  Let us now turn to the latter inequality in~\eqref{eq:vk-bounded}.  It follows from the normal direction computation that $\|c_k\| \geq \|c_k+J_k v_k\|$, which by the triangle inequality implies that $\|J_kv_k\| \leq 2\|c_k\|$.  Note that since $v_k\in\Range(J_k^T)$, one has $v_k = J_k^Tw_k$ where $w_k = (J_kJ_k^T)^{-1}J_kv_k$.  Putting these facts together shows that
  \bequationNN
    \|v_k\| = \|J_k^T w_k\| = \|J_k^T(J_kJ_k^T)^{-1}J_kv_k\| \leq \|J_k^T\|\|(J_kJ_k^T)^{-1}\|\|J_kv_k\| \leq \tfrac{2\|J_k^T\|}{\lambdamin}\|c_k\|,
  \eequationNN
  which combined with Assumption~\ref{ass.main} establishes the existence of a $\omega_3\in\R{}_{>0}$ such that the second inequality in~\eqref{eq:vk-bounded} holds, as desired.
\eproof

The next result gives a bound on the size of the search direction relative to the constraint violation and the size of the normal step.

\blemma\label{lm:Thetak_lower_bound}
  There exists $\omega_4\in\R{}_{\geq 2}$ such that, for all $k\in\N{}$, it holds that
  \bequationNN
    \|d_k\|^2 \leq \omega_4(\|u_k\|^2 + \|c_k\|).
  \eequationNN
\elemma
\bproof
  Observe that $0 \leq (\|u_k\|-\|v_k\|)^2 = \|u_k\|^2 + \|v_k\|^2 - 2\|u_k\|\|v_k\|$.  Using this fact, $d_k = v_k + u_k$, the triangle inequality, and Lemma~\ref{lm:vk_bounded}, it follows that
  \bequationNN
    \baligned
      \|d_k\|^2 
        &\leq (\|u_k\| + \|v_k\|)^2 = \|u_k\|^2 + \|v_k\|^2 + 2\|u_k\|\|v_k\| \\
        &\leq 2(\|u_k\|^2 + \|v_k\|^2) \leq 2(\|u_k\|^2 + \omega_3^2\|c_k\|^2) \\
        &\leq \max\{2,2\omega_3^2\|c_k\|\}(\|u_k\|^2 + \|c_k\|).
    \ealigned
  \eequationNN
  The existence of the required $\omega_4\in\R{}_{\geq 2}$ now follows from Assumption~\ref{ass.main} since $\max\{2,2\omega_3^2\|c_k\|\}$ is uniformly bounded for all $k\in\N{}$, which completes the proof.
\eproof

The next lemma shows that the model reduction $\Delta l(x_k,\tau_k,g_k,v_k+u_k)$ is bounded below by a similar quantity as the upper bound for $\|d_k\|^2$ in the previous lemma.

\blemma\label{lm:Thetak_upper_bound}
  There exists $\kappa_l\in\R{}_{>0}$ such that for all $k\in\N{}$, it holds that
  \bequationNN
    \Delta l(x_k,\tau_k,g_k,v_k + u_k) \geq \kappa_l\tau_k(\|u_k\|^2 + \|c_k\|) \geq \tfrac{\kappa_l\tau_k}{\omega_4}\|d_k\|^2 > 0.
  \eequationNN
\elemma
\bproof
  Lemma~\ref{lm:model_reduction_true} shows that \eqref{eq.model_reduction_condition} holds with $\tau \equiv \tau_k$.  Combining this fact with Lemma~\ref{lm:feasibility_decrease} and the monotonically nonincreasing behavior of $\{\tau_k\}$ shows that
  \bequationNN
    \baligned
      \Delta l(x_k,\tau_k,g_k,v_k + u_k) &\geq \sigma_u\tau_k\max\{u_k^TH_ku_k,\epsilon_u\|u_k\|^2\} + \sigma_c (\|c_k\| - \|c_k + J_kv_k\|) \\
      &\geq \sigma_u\tau_k\epsilon_u\|u_k\|^2 + \sigma_c\omega_1\|c_k\| \geq \tau_k( \sigma_u\epsilon_u\|u_k\|^2 + \sigma_c\omega_1\|c_k\|/\tau_{-1}) \\
      &\geq \min\{\sigma_u\epsilon_u,\tfrac{\sigma_c\omega_1}{\tau_{-1}}\} \tau_k ( \|u_k\|^2 + \|c_k\|),
    \ealigned
  \eequationNN
  which proves the existence of the claimed $\kappa_l \in\R{}_{>0}$ since $\sigma_u$, $\epsilon_u$, $\sigma_c$, $\omega_1$, and $\tau_{-1}$ are positive real numbers.  The remaining inequalites follow from Lemmas~\ref{lm:Thetak_lower_bound} and \ref{lem:d-nonzero}.
\eproof

We next prove a deterministic uniform lower bound for the sequence $\{\xi_k\}$.

\blemma\label{lem:xi-bound}
  There exists $\xi_{\min}\in\R{}_{>0}$ such that, in any run of the algorithm, there exists $k_\xi\in\N{}$ and $\xi_{k_\xi}\in[\xi_{\min},\infty)$ such that $\xi_k = \xi_{k_\xi}$ for all $k\geq k_\xi$.
\elemma
\bproof
  For all $k\in\N{}$, it follows from~\eqref{eq:stochastic_ratio} and Lemmas~\ref{lm:Thetak_lower_bound} and \ref{lm:Thetak_upper_bound} that
  \bequation\label{xitrial-bd}
    \xitrial_k = \tfrac{\Delta l(x_k,\tau_k,g_k,d_k)}{\tau_k\|d_k\|^2} \geq \tfrac{\kappa_l\tau_k(\|u_k\|^2 + \|c_k\|)}{\tau_k\omega_4(\|u_k\|^2 + \|c_k\|)} = \tfrac{\kappa_l}{\omega_4}. 
  \eequation
  Now, consider any iteration such that $\xi_k < \xi_{k-1}$.  For such iterations, it follows from~\eqref{eq:stochastic_ratio} and~\eqref{xitrial-bd} that $\xi_k \geq (1-\epsilon_\xi)\xitrial_k \geq (1-\epsilon_{\xi})\kappa_l/\omega_4$.  Combining this fact with the initial choice of $\xi_{-1}$ shows that $\xi_k \geq \xi_{\min} := \min\{(1-\epsilon_{\xi})\kappa_l/\omega_4, \xi_{-1}\}$ for all $k\in\N{}$. Combining this result with the fact that anytime $\xi_k < \xi_{k-1}$ it must hold that $\xi_k \leq (1-\epsilon_\xi)\xi_{k-1}$ (it decreases by at least a factor of $1-\epsilon_\xi$), gives the desired result.
\eproof

The next lemma gives a bound on the change in the merit function each iteration.

\blemma\label{lm:stochastic_merit_decrease_upper_bound}
  For all $k\in\N{}$, it holds that
  \bequationNN
    \baligned
      &\ \ \ \ \phi(x_k + \alpha_kd_k,\tau_k) - \phi(x_k,\tau_k) \\
      &\leq -\alpha_k\Delta l(x_k,\tau_k,\nabla f_k,\dTrue_k) + \alpha_k\tau_k\nabla f_k^T(d_k - \dTrue_k)
+ (1-\eta)\alpha_k\beta_k\Delta l(x_k,\tau_k,g_k,d_k) \\
      &\ \ \ + \alpha_k\|c_k + J_kd_k\| - \alpha_k\|c_k + J_kv_k\|.
    \ealigned
  \eequationNN
\elemma
\bproof
  By Lemma~\ref{lm:stepsize_well_defined}, one has that $\varphi(\alpha_k) \leq 0$.  Hence, starting with the third row of~\eqref{eq.merit_reduction_bound}, adding and subtracting the terms $\alpha_k\tau_k\nabla f_k^T\dTrue_k$, $\alpha_k \|c_k\|$, $\alpha_k \|c_k + J_k\dTrue_k\|$, and $\alpha_k\beta_k\Delta l(x_k,\tau_k,g_k,d_k)$, using the definition of $\varphi(\cdot)$, and using the fact that $J_k \dTrue_k = J_k (v_k + \uTrue_k) = J_k v_k$, one finds that
  \bequationNN
    \baligned
      &\ \phi(x + \alpha_kd_k, \tau_k) - \phi(x_k,\tau_k) \\
  \leq&\ \alpha_k\tau_k\nabla f_k^Td_k + \|c_k+\alpha_k J_kd_k\| - \|c_k\| + \tfrac{1}{2}(\tau_kL+\Gamma)\alpha_k^2\|d_k\|^2 \\
     =&\ -\alpha_k\Delta l(x_k,\tau_k,\nabla f_k,\dTrue_k) + \alpha_k\tau_k\nabla f_k^T(d_k - \dTrue_k) + (\alpha_k - 1)\|c_k\| \\
      &\ + \|c_k+\alpha_kJ_kd_k\| - \alpha_k\|c_k + J_k\dTrue_k\| + \tfrac{1}{2}(\tau_kL+\Gamma)\alpha_k^2\|d_k\|^2 \\
      &\ - \alpha_k\beta_k\Delta l(x_k,\tau_k,g_k,d_k) + \alpha_k\beta_k\Delta l(x_k,\tau_k,g_k,d_k) \\
  \leq&\ -\alpha_k\Delta l(x_k,\tau_k,\nabla f_k,\dTrue_k) + \alpha_k\tau_k\nabla f_k^T(d_k - \dTrue_k) + \alpha_k\|c_k + J_kd_k\| \\
      &\ - \alpha_k\|c_k + J_k\dTrue_k\| - \eta\alpha_k\beta_k\Delta l(x_k,\tau_k,g_k,d_k) + \alpha_k\beta_k\Delta l(x_k,\tau_k,g_k,d_k) \\
     =&\ -\alpha_k\Delta l(x_k,\tau_k,\nabla f_k,\dTrue_k) + \alpha_k\tau_k\nabla f_k^T(d_k - \dTrue_k) + (1-\eta)\alpha_k\beta_k\Delta l(x_k,\tau_k,g_k,d_k) \\
      &\ + \alpha_k\|c_k + J_kd_k\| - \alpha_k\|c_k + J_kv_k\|, \\
    \ealigned
  \eequationNN
  which completes the proof. 
\eproof

We now derive bounds on the expected difference between $u_k$ and $\uTrue_k$.  To that end, let us define $Z_k\in\R{n\times (n-m)}$ as a matrix whose columns form an orthonormal basis for $\Null(J_k)$, which implies that $Z_k^T Z_k = I$ and $J_k Z_k = 0$.  Under Assumption~\ref{ass.main}, let $u_{k,1}\in\R{m}$ and $u_{k,2}\in\R{n-m}$ be vectors forming the orthogonal decomposition of $u_k$ into $\Range(J_k^T)$ and $\Null(J_k)$ in the sense that $u_k = J_k^Tu_{k,1} + Z_ku_{k,2}$.  It follows from~\eqref{eq.linear_system} that $u_{k,1} = (J_kJ_k^T)^{-1}r_k$ and $u_{k,2} = -(Z_k^TH_kZ_k)^{-1}Z_k^T(g_k + H_kv_k + H_kJ_k^T(J_kJ_k^T)^{-1}r_k - \rho_k)$, with which one can derive:
\bequation\label{eq:u_solution}
  \baligned
    u_k &= J_k^T(J_kJ_k^T)^{-1}r_k - Z_k(Z_k^TH_kZ_k)^{-1}Z_k^T(g_k + H_kv_k + H_kJ_k^T(J_kJ_k^T)^{-1}r_k - \rho_k) \\
    \uTrue_k &= -Z_k(Z_k^TH_kZ_k)^{-1}Z_k^T(\nabla f_k + H_kv_k).
  \ealigned
\eequation
The corresponding values for $\deltaTrue_k$ and $\delta_k$ are found to be:
\bequation\label{eq:delta_solution}
  \baligned
    \delta_k &= -(J_kJ_k^T)^{-1}J_k(g_k + H_kv_k + H_ku_k - \rho_k) - y_k \\
    \deltaTrue_k &= -(J_kJ_k^T)^{-1}J_k(\nabla f_k + H_kv_k + H_k\uTrue_k) - y_k. 
  \ealigned
\eequation
In the proof of the lemma below, we use the fact that
\bequation\label{Z-bound}
  \|I - Z_k(Z_k^TH_kZ_k)^{-1}Z_k^TH_k\| \leq 1,
\eequation
which can be seen as follows: The nonzero eigenvalues of a matrix product $AB$ are equal to the nonzero eigenvalues of $BA$ when the product is valid, from which it follows that the nonzero eigenvalues of $Z_k(Z_k^TH_kZ_k)^{-1}Z_k^TH_k$ are precisely the eigenvalues of $Z_k^TH_k Z_k(Z_k^TH_kZ_k)^{-1} = I$, which are all equal to one; hence, the bound in~\eqref{Z-bound} holds.

\blemma\label{lem:u-diff}
  There exists $\omega_5 \in \R{}_{>0}$ such that, for all $k \in \N{}$, it holds that
  \bequationNN
    \|\E_k[u_k - \uTrue_k]\| \leq \omega_5\beta_k \ \ \text{and}\ \ \E_k[\|u_k - \uTrue_k\|] \leq \zeta^{-1}\sqrt{M_g} + \omega_5\beta_k
  \eequationNN
\elemma
\bproof
  It follows from~\eqref{eq:u_solution} that
  \bequationNN
    u_k - \uTrue_k = J_k^T(J_kJ_k^T)^{-1}r_k - Z_k(Z_k^TH_kZ_k)^{-1}Z_k^T(g_k - \nabla f_k + H_kJ_k^T(J_kJ_k^T)^{-1}r_k - \rho_k),
  \eequationNN
  which combined with Assumption~\ref{ass.g} shows that
  \bequationNN
    \baligned
      \E_k[u_k - \uTrue_k] = &(I - Z_k(Z_k^TH_kZ_k)^{-1} Z_k^TH_k ) J_k^T(J_kJ_k^T)^{-1} \E_k[r_k] \\
      &+ Z_k(Z_k^TH_kZ_k)^{-1}Z_k^T\E_k[\rho_k].
    \ealigned
  \eequationNN
  Combining this equation with the triangle inequality, Assumptions~\ref{ass.H} and~\ref{ass.main}, \eqref{eq.pd_residual_condition}, and~\eqref{Z-bound} ensures the existence of $\omega_5\in\R{}_{>0}$ such that, for all $k \in \N{}$,
  \bequationNN
    \baligned
      \|\E_k[ u_k - \uTrue_k ]\| 
        &\leq \|J_k^T(J_kJ_k^T)^{-1}\|\|\E_k[r_k]\| + \zeta^{-1}\|\E_k[\rho_k]\| \\
        &\leq \|J_k^T(J_kJ_k^T)^{-1}\|\kappa_r\beta_k + \zeta^{-1}\kappa_{\rho}\beta_k \leq \omega_5\beta_k,
    \ealigned
  \eequationNN
  is the first desired result.  Next, to derive the desired bound on $\E_k[\|u_k - \uTrue_k\|]$, one can combine the expression above for $u_k - \uTrue_k$ with the triangle inequality to obtain 
  \bequationNN
    \baligned
      \|u_k - \uTrue_k\| \leq&\ \|Z_k(Z_k^TH_kZ_k)^{-1}Z_k^T(g_k - \nabla f_k)\| + \|Z_k(Z_k^TH_kZ_k)^{-1}Z_k^T\rho_k\| \\
      &+ \|(I - Z_k(Z_k^TH_kZ_k)^{-1}Z_k^TH_k)J_k^T(J_kJ_k^T)^{-1}r_k\|.
    \ealigned
  \eequationNN
  Taking conditional expectation and using Assumption~\ref{ass.g}, \eqref{jensens-M}, \eqref{Z-bound}, and~\eqref{eq.pd_residual_condition},
  \bequationNN
    \baligned
      \E_k[\|u_k - \uTrue_k\|] &\leq \zeta^{-1}\sqrt{M_g} + \zeta^{-1}\E_k[\|\rho_k\|] + \|J_k^T(J_kJ_k^T)^{-1}\|\E_k[\|r_k\|] \\
      &\leq \zeta^{-1}\sqrt{M_g} + \zeta^{-1}\kappa_{\rho}\beta_k + \|J_k^T(J_kJ_k^T)^{-1}\|\kappa_r\beta_k \leq \zeta^{-1}\sqrt{M_g} + \omega_5\beta_k,
    \ealigned
  \eequationNN
  where $\omega_5$ is the same value as used above, which completes the proof.
\eproof

We now bound the difference (in expectation) between $\nabla f_k^T\dTrue_k$ and $g_k^Td_k$.

\blemma\label{lem:gd_diff}
  There exist $(\omega_6,\omega_7) \in \R{}_{>0} \times \R{}_{>0}$ such that, for all $k\in\N{}$,
  \bequationNN
    |\E_k[\nabla f_k^T\dTrue_k - g_k^Td_k]| \leq \omega_6\beta_k + \omega_7\beta_k\sqrt{M_g} + \zeta^{-1}M_g.
  \eequationNN
\elemma
\bproof
  It follows from the triangle inequality and linearity of $E_k$ that
  \bequationNN
    \baligned
      |\E_k[\nabla f_k^T\dTrue_k - g_k^Td_k]|
        &= |\E_k[\nabla f_k^T(\dTrue_k - d_k) + (\nabla f_k - g_k)^Td_k]| \\
        &\leq |\nabla f_k^T\E_k[\dTrue_k - d_k]| + |\E_k[(\nabla f_k - g_k)^Td_k]|.
    \ealigned
  \eequationNN
  For the first term on the right-hand side, it follows by the Cauchy-Schwarz inequality, $\dTrue_k = v_k + \uTrue_k$, $d_k = v_k + u_k$, and Lemma~\ref{lem:u-diff} that there exists $\omega_6 \in\R{}_{>0}$ with
  \bequationNN
    \baligned
      |\nabla f_k^T\E_k[\dTrue_k - d_k]| &\leq \|\nabla f_k\|\|\E_k[\dTrue_k - d_k]\| \\
        &= \|\nabla f_k\|\|\E_k[\uTrue_k - u_k]\| \leq \omega_6\beta_k. 
    \ealigned
  \eequationNN
  For the second term on the right-hand side, first observe from Assumption~\ref{ass.g} that $\E_k[(\nabla f_k-g_k)^T v_k] = v_k^T \E_k[\nabla f_k-g_k] = 0$.  Combining this fact with~\eqref{eq:u_solution}, the triangle inequality, the Cauchy-Schwarz inequality, Assumptions~\ref{ass.main}--\ref{ass.H}, \eqref{Z-bound}, and~\eqref{jensens-M} shows that there exist $(\bar\omega_7,\omega_7) \in \R{}_{>0} \times \R{}_{>0}$ such that
  \bequationNN
    \baligned
       &\ |\E_k[(\nabla f_k - g_k)^Td_k]| \\
      =&\ |\E_k[(\nabla f_k - g_k)^T((I - Z_k(Z_k^TH_kZ_k)^{-1}Z_k^TH_k)J_k^T(J_kJ_k^T)^{-1}r_k \\
       &\ -Z_k(Z_k^TH_kZ_k)^{-1}Z_k^T(g_k - \nabla f_k - \rho_k))]| \\
   \leq&\ |\E_k[(\nabla f_k - g_k)^T(I - Z_k(Z_k^TH_kZ_k)^{-1}Z_k^TH_k)J_k^T(J_kJ_k^T)^{-1}r_k]| \\ 
       &\ + |\E_k[(\nabla f_k - g_k)^TZ_k(Z_k^TH_kZ_k)^{-1}Z_k^T\rho_k]| \\
       &\ + |\E_k[(\nabla f_k - g_k)^TZ_k(Z_k^TH_kZ_k)^{-1}Z_k^T(\nabla f_k - g_k)]| \\
   \leq&\ \E_k[\|\nabla f_k - g_k\|\|(I - Z_k(Z_k^TH_kZ_k)^{-1}Z_k^TH_k)J_k^T(J_kJ_k^T)^{-1}\| \|r_k\|] \\
       &\ + \E_k[\|\nabla f_k - g_k\| \|Z_k(Z_k^TH_kZ_k)^{-1}Z_k^T\| \|\rho_k\|] + \zeta^{-1}\E_k[\|\nabla f_k - g_k\|^2] \\
   \leq&\ (\bar\omega_7\kappa_r\beta_k + \zeta^{-1}\kappa_{\rho}\beta_k)\E_k[\|\nabla f_k - g_k\|] + \zeta^{-1}M_g \\
   \leq&\ (\bar\omega_7\kappa_r + \zeta^{-1}\kappa_{\rho})\beta_k\sqrt{M_g} + \zeta^{-1}M_g = \omega_7\beta_k\sqrt{M_g} + \zeta^{-1}M_g.
    \ealigned
  \eequationNN
  Combining the results above gives the desired result.
\eproof

We now proceed to bound (in expectation) the last two terms appearing in the right-hand side of the inequality proved in Lemma~\ref{lm:stochastic_merit_decrease_upper_bound}.

\blemma\label{lem:Jd_Jv_diff}
  There exists $\omega_8 \in \R{}_{>0}$ such that, for all $k\in\N{}$, it holds that
  \bequationNN
    \E_k[\alpha_k(\|c_k + J_kd_k\| - \|c_k + J_kv_k\|)] \leq \omega_8\beta_k^2.
  \eequationNN
\elemma
\bproof
  From the triangle inequality, \eqref{eq.linear_system}, \eqref{eq.pd_residual_condition}, the fact that $\alpha_k \in [\alphamin_k,\alphamax_k]$, \eqref{def-alphamax}, \eqref{eq:stochastic_alpha_suff.new}, \eqref{betak-requirement}, and the monotonically nonincreasing behavior of $\{\tau_k\}$ and $\{\xi_k\}$, it follows that there exists $\omega_8 \in \R{}_{>0}$ such that
  \bequationNN
    \baligned
      &\ \E_k[\alpha_k(\|c_k + J_kd_k\| - \|c_k + J_kv_k\|)] \leq \E_k[\alpha_k\|J_ku_k\|] = \E_k[\alpha_k\|r_k\|] \\
  \leq&\ \kappa_r\beta_k\E_k[\alphamax_k] \leq \kappa_r\beta_k\E_k[\alphamin_k + \theta\beta_k^2] = \kappa_r\beta_k \E_k[ (\tfrac{2(1-\eta)\beta_k\xi_k\tau_k}{\tau_kL+\Gamma} + \theta\beta_k^2)] \\
  \leq&\ (\tfrac{2(1-\eta)\xi_{-1}\tau_{-1}}{\Gamma} + \theta\beta_k)\kappa_r\beta_k^2 \leq \omega_8\beta_k^2,
    \ealigned
  \eequationNN
  which gives the desired conclusion.
\eproof

\subsection{Convergence analysis}\label{sec.convergence}

Our goal now is to prove a convergence result for our algorithm.  In general, in a run of the algorithm, one of three possible events can occur.  One possible event is that the merit parameter sequence eventually remains constant at a value that is \emph{sufficiently small}.  This is the event that we consider in our analysis here, where the meaning of \emph{sufficiently small} is defined formally below.  The other two possible events are that the merit parameter sequence vanishes or eventually remains constant at a value that is too large.  As discussed in \cite[Section~3.2.2]{BeraCurtRobiZhou21}, the former of these two events does not occur if the differences between the stochastic gradient estimates and the true gradients of the objective remain uniformly bounded in norm, and the latter of these two events occurs with probability zero in a given run of the algorithm if one makes a reasonable assumption about the influence of the stochastic gradient estimates on the computed search directions; see also \cite[Section~4.3]{BeraCurtONeiRobi21} for additional discussion of the latter case in the context of an algorithm that employs a step decomposition approach, as does our algorithm.  For our purposes here, we do not consider these latter two events since we contend that, for practical purposes, they can be ignored for the same reasons as are claimed in \cite{BeraCurtRobiZhou21}.

To define our event of interest, consider for each $k\in\N{}$ the condition
\bequation\label{check-truetrial}
  \nabla f_k^T\dTrue_k + \max\{(\uTrue_k)^TH_k\uTrue_k,\epsilon_u \|\uTrue_k\|^2\} \leq 0
\eequation
(similar to the one appearing in~\eqref{tautrail-def}).  With this condition, let us define the following trial value of the merit parameter that would be computed in iteration $k \in \N{}$ (conditioned on $x_k$ being the $k$th iterate) if the algorithm were to employ $\nabla f_k$ in place of $g_k$ and compute an exact solution of the linear system \eqref{prob.tangential_true}:
\bequationNN
  \tauTruetrial_k \gets
  \bcases
    \infty & \text{if \eqref{check-truetrial} holds,} \\
    \frac{(1-\tfrac{\sigma_c}{\epsilon_r})(\|c_k\|-\|c_k + J_k\dTrue_k\|)}{\nabla f_k^T\dTrue_k + \max\{(\uTrue_k)^TH_k\uTrue_k,\epsilon_u \|\uTrue_k\|^2\}} & \text{if \eqref{check-truetrial} does not hold.}
  \ecases
\eequationNN
(To be clear, the quantity $\tauTruetrial_k$ never needs to be computed by our algorithm; it is only used in our analysis in this subsection.)  Using this quantity, we define our event of interest, namely, $E_{\tau,\low}$, as the following.

\begin{tcolorbox}[colback=white]
\textbf{Event $E_{\tau,\low}$.}
  Event $E_{\tau,\low}$ occurs if and only if there exists an iteration number $k_{\tau,\xi} \in \N{}$ such that, with $\xi_{\min}$ given in Lemma~\ref{lem:xi-bound}, it holds that
  \bequation\label{eq:stochastic_good_case}
    \tau_k = \tau_{k_{\tau,\xi}} \leq \tauTruetrial_k\ \ \text{and}\ \ \xi_k = \xi_{k_{\tau,\xi}} \geq \xi_{\min}\ \ \text{for all}\ \ k\geq k_{\tau,\xi}.
  \eequation
\end{tcolorbox}  

For our analysis in this subsection, the following supersedes Assumption~\ref{ass.g}.

\bassumption\label{ass.tau_low}
  There exists $M_g \in \R{}_{>0}$ such that, for all $k \in \N{}$, the stochastic gradient $g_k$ has the properties that $\E_{k,\tau,\low}[g_k] = \nabla f_k$ and $\E_{k,\tau,\low}[\|g_k - \nabla f_k\|_2^2] \leq M_g$, where $\E_{k,\tau,\low}[\cdot]$ denotes expectation with respect to the distribution of $\omega$ conditioned on the event that $E_{\tau,\low}$ occurs and $x_k$ is the primal iterate in iteration $k \in \N{}$.
\eassumption

Our results in this subsection focus on $k \in \N{}$ with $k \geq k_{\tau,\xi}+1$, at which point, in any run in which Event~$E_{\tau,\low}$ occurs, the merit parameter satisfies $\tau_k = \tau_{k_{\tau,\xi}}$ independently from the stochastic gradient $g_k$ that is generated.

Our first result provides an upper bound (in expectation) for the second term appearing on the right-hand side of the inequality in Lemma~\ref{lm:stochastic_merit_decrease_upper_bound}.

\blemma\label{lm:stochastic_gd_diff}
  Under Event~$E_{\tau,\low}$, there exists $\omega_9\in\R{}_{>0}$ such that
  \bequationNN
    \E_{k,\tau,\low}[\alpha_k\tau_k\nabla f_k^T(d_k - \dTrue_k)] \leq \omega_9\beta_k^2\ \ \text{for all}\ \ k \geq k_{\tau,\xi} + 1.
  \eequationNN
\elemma
\bproof
  Under Assumption~\ref{ass.tau_low}, the logic as in the proof of Lemma~\ref{lem:u-diff} allows us to conclude that, under $E_{\tau,low}$, it holds for all $k \in \N{}$ that
  \bequation\label{new-E}
    \|\E_{k,\tau,\low}[u_k - \uTrue_k]\| \leq \omega_5\beta_k \ \ \text{and}\ \ \E_{k,\tau,\low}[\|u_k - \uTrue_k\|] \leq \zeta^{-1}\sqrt{M_g} + \omega_5\beta_k.
  \eequation
  Let $E_k$ be the event that $\nabla f_k^T(d_k - \dTrue_k) \geq 0$ and let $E_k^c$ be its complementary event. Let $\Pmbb_{k,\tau,\low}[\cdot]$ denote probability conditioned on the occurrence of event $E_{\tau,\low}$ and $x_k$ being the $k$th primal iterate.  It now follows from~\eqref{eq:stochastic_good_case}, the definition of $E_k$, the fact that $\alpha_k \in[\alphamin_k,\alphamax_k]$, and the Law of Total Expectation that for all $k \geq k_{\tau,\xi}+1$
  \bequationNN
    \baligned
      &\ \E_{k,\tau,\low}[\alpha_k\tau_k\nabla f_k^T(d_k-\dTrue_k)] \\
     =&\ \E_{k,\tau,\low}[\alpha_k\tau_{k_{\tau,\xi}}\nabla f_k^T(d_k-\dTrue_k) | E_k]\Pmbb_{k,\tau,\low}[E_k] \\
      &\ + \E_{k,\tau,\low}[\alpha_k\tau_{k_{\tau,\xi}}\nabla f_k^T(d_k-\dTrue_k) | E_k^c]\Pmbb_{k,\tau,\low}[E_k^c] \\
  \leq&\ \E_{k,\tau,\low}[\alphamax_k\tau_{k_{\tau,\xi}}\nabla f_k^T(d_k-\dTrue_k) | E_k]\Pmbb_{k,\tau,\low}[E_k] \\
      &\ + \E_{k,\tau,\low}[\alphamin_k\tau_{k_{\tau,\xi}}\nabla f_k^T(d_k-\dTrue_k) | E_k^c]\Pmbb_{k,\tau,\low}[E_k^c] \\
     =&\ \E_{k,\tau,\low}[(\alphamax_k - \alphamin_k)\tau_{k_{\tau,\xi}}\nabla f_k^T(d_k-\dTrue_k) | E_k]\Pmbb_{k,\tau,\low}[E_k] \\
      &\ + \E_{k,\tau,\low}[\alphamin_k\tau_{k_{\tau,\xi}}\nabla f_k^T(d_k-\dTrue_k)].
    \ealigned
  \eequationNN
  Combining this with the fact that \eqref{def-alphamax} ensures $\alphamax_k-\alphamin_k \leq \theta \beta_k^2$, the Cauchy-Schwarz inequality, the fact that $\alphamin_k = 2(1-\eta)\beta_k\xi_{k_{\tau,\xi}}\tau_{k_{\tau,\xi}}/(\tau_{k_{\tau,\xi}}L+\Gamma)$ for all $k \geq k_{\tau,\xi} + 1$, and the Law of Total Expectation shows for all $k \geq k_{\tau,\xi}+1$ that
  \bequationNN
    \baligned
      \E_{k,\tau,\low}[\alpha_k\tau_k\nabla f_k^T(d_k-\dTrue_k)]
        \leq&\ \theta\beta_k^2\tau_{k_{\tau,\xi}}\|\nabla f_k\| \E_{k,\tau,\low}[\|d_k-\dTrue_k\| | E_k]\Pmbb_{k,\tau,\low}[E_k] \\
            &\ + \tfrac{2(1-\eta)\beta_k\xi_{k_{\tau,\xi}}\tau_{k_{\tau,\xi}}}{\tau_{k_{\tau,\xi}}L+\Gamma}
\tau_{k_{\tau,\xi}}\|\nabla f_k\| \| \E_{k,\tau,\low}[d_k - \dTrue_k]\| \\
        \leq&\ \theta\beta_k^2\tau_{k_{\tau,\xi}}\|\nabla f_k\| \E_{k,\tau,\low}[\|d_k-\dTrue_k\|] \\
            &\ + \tfrac{2(1-\eta)\beta_k\xi_{k_{\tau,\xi}}\tau_{k_{\tau,\xi}}}{\tau_{k_{\tau,\xi}}L+\Gamma} \tau_{k_{\tau,\xi}}\|\nabla f_k\| \| \E_{k,\tau,\low}[d_k - \dTrue_k]\|
    \ealigned
  \eequationNN
  Combining this with~\eqref{new-E}, \eqref{betak-requirement}, $\|d_k - \dTrue_k\| = \|v_k+u_k-(v_k+\uTrue_k)\| = \|u_k-\uTrue_k\|$, and Assumption~\ref{ass.main} shows there exists $\omega_9 \in \R{}_{>0}$ where, for all $k \geq k_{\tau,\xi}+1$,
  \bequationNN
    \baligned
      &\ \E_{k,\tau,\low}[\alpha_k\tau_k\nabla f_k^T(d_k-\dTrue_k)] \\
  \leq&\ \theta\beta_k^2\tau_{k_{\tau,\xi}}\|\nabla f_k\| (\zeta^{-1}\sqrt{M_g} + \omega_5\beta_k)  
+ \tfrac{2(1-\eta)\beta_k\xi_{k_\xi}\tau_{k_{\tau,\xi}}}{\tau_{k_{\tau,\xi}}L+\Gamma}
\tau_{k_{\tau,\xi}}\|\nabla f_k\| \omega_5\beta_k 
\leq \omega_9\beta_k^2,
    \ealigned
  \eequationNN
  which is the desired conclusion.
\eproof

We now use the model reduction based on the true step $\dTrue_k$ to build an upper bound on the (expected) reduction in the model based on the step $d_k$.

\blemma\label{lm:stochastic_Delta_q}
  Under Event~$E_{\tau,\low}$, it holds for all $k\geq k_{\tau,\xi} + 1$ that
  \bequationNN
    \baligned
      &\ \E_{k,\tau,\low}[\Delta l(x_k,\tau_k,g_k,d_k)] \\
  \leq&\ \Delta l(x_k,\tau_{k_{\tau,\xi}},\nabla f_k,\dTrue_k) + \kappa_r\beta_k + \tau_{k_{\tau,\xi}}(\omega_6\beta_k + \omega_7\beta_k\sqrt{M_g} + \zeta^{-1}M_g).
    \ealigned
  \eequationNN
\elemma
\bproof
  Under Assumption~\ref{ass.tau_low}, the logic as in the proof of Lemma~\ref{lem:gd_diff} allows us to conclude that, under $E_{\tau,\low}$, it holds for all $k \in \N{}$ that
  \bequationNN
    |\E_{k,\tau,\low}[\nabla f_k^T\dTrue_k - g_k^Td_k]| \leq \omega_6\beta_k + \omega_7\beta_k\sqrt{M_g} + \zeta^{-1}M_g.
  \eequationNN
  It follows from this, \eqref{eq.model_reduction}, the fact that $d_k = v_k + u_k$, the triangle inequality, the fact that $c_k$, $J_k$, $v_k$, $\nabla f_k$, and $\dTrue_k$ are all deterministic conditioned on $x_k$ as the $k$th primal iterate, \eqref{eq.linear_system_true}, and \eqref{eq.pd_residual_condition} that for all $k\geq k_{\tau,\xi} + 1$
  \bequationNN
    \baligned
      &\ \E_{k,\tau,\low}[\Delta l(x_k,\tau_k,g_k,d_k)] = \E_{k,\tau,\low}[-\tau_{k_{\tau,\xi}}g_k^Td_k + \|c_k\| - \|c_k + J_kd_k\|] \\
  \leq&\ \Delta l(x_k,\tau_{k_{\tau,\xi}},\nabla f_k,\dTrue_k) + \kappa_r\beta_k + \tau_{k_{\tau,\xi}}(\omega_6\beta_k + \omega_7\beta_k\sqrt{M_g} + \zeta^{-1}M_g),
    \ealigned
  \eequationNN
  which is the desired result.
\eproof

For the final result of this section, we define
\bequation\label{def:Etaulow}
  \E_{\tau,\low}[\cdot] = \E[\ \cdot \ | \ \text{Event $E_{\tau,\low}$ occurs and Assumption~\ref{ass.tau_low} holds}].
\eequation
In the result, the quantity $\Delta l(x_k,\tau_k,\nabla f_k, \dTrue_k)$ serves as a measure of stationarity with respect to \eqref{prob.opt}; after all, the proof for Lemma~\ref{lm:Thetak_upper_bound} shows, with $(\nabla f_k,\uTrue_k,\dTrue_k)$ in place of $(g_k,u_k,d_k)$, that by Assumption~\ref{ass.tau_low} it follows for  $k\geq k_{\tau,\xi}+1$ that
\bequation\label{Deltal-true}
  \Delta l(x_k,\tau_{k_{\tau,\xi}},\nabla f_k,\dTrue_k) \geq \kappa_l\tau_{k_{\tau,\xi}}(\|\uTrue_k\|^2 + \|c_k\|) \geq \tfrac{\kappa_l\tau_{k_{\tau,\xi}}}{\omega_4}\|\dTrue_k\|^2 > 0.
\eequation
Thus, if there is an infinite $\Kcal \subseteq \N{}$ with $\lim_{k \in \Kcal, k \to \infty} \Delta l(x_k,\tau_{k_{\tau,\xi}},\nabla f_k,\dTrue_k) = 0$, then it follows from~\eqref{Deltal-true} and~\eqref{lm:vk_bounded} that $\lim_{k\in\Kcal,k\to\infty} \|c_k\| =  \lim_{k\in\Kcal,k\to\infty} \|\uTrue_k\| = \lim_{k\in\Kcal,k\to\infty} \|v_k\| = 0$, which combined with~\eqref{eq.linear_system_true} shows that any limit point of $\{(x_k,y_k+\deltaTrue_k)\}$ is a first-order stationary point for~\eqref{prob.opt}.  In our stochastic setting, we cannot guarantee that such a limit holds surely.  Rather, in the following result, we prove for two different choices of $\{\beta_k\}$ that an expected average of this measure of stationarity exhibits desirable properties.  These properties match those ensured by a stochastic gradient method in the unconstrained setting (where $\|\nabla f_k\|^2$ plays the role of the measure of stationarity for the minimization of $f$).

\btheorem\label{thm:stochastic_good_case_final}
  Under Event~$E_{\tau,\low}$, let $k_{\tau,\xi}$ be defined as in~\eqref{eq:stochastic_good_case} and define $\Abar = \tfrac{2(1-\eta)\xi_{\min}\tau_{k_{\tau,\xi}}}{\tau_{k_{\tau,\xi}}L+\Gamma}$ and $\Mbar = (1-\eta)(\Abar+\theta)\big(\kappa_r + \tau_{k_{\tau,\xi}}(\omega_6 + \omega_7\sqrt{M_g} + \zeta^{-1}M_g)\big) + \omega_8 + \omega_9$, where $\xi_{\min}$ is defined in Lemma~\ref{lem:xi-bound}. Then, the following results hold:
  \bitemize
    \item[(i)] If $\beta_k = \beta \in (0,\Abar/((1-\eta)(\Abar + \theta)))$ for all $k\geq k_{\tau,\xi} + 1$, then
    \bequation\label{eq:main_thm_fixed_stepsize}
      \baligned
        &\ \E_{\tau,\low}\left[\frac{1}{K}\sum_{j=k_{\tau,\xi}+1}^{k_{\tau,\xi}+K} \Delta l(x_j,\tau_{k_{\tau,\xi}},\nabla f_j,\dTrue_j)\right] \\
    \leq&\ \tfrac{\beta\Mbar}{\Abar - (1-\eta)(\Abar + \theta)\beta} + \tfrac{\E_{\tau,\low}[\phi(x_{k_{\tau,\xi}+1},\tau_{k_{\tau,\xi}})] - \phi_{\min}}{K\beta(\Abar - (1-\eta)(\Abar + \theta)\beta)} \xrightarrow{K\to\infty} \tfrac{\beta\Mbar}{\Abar - (1-\eta)(\Abar + \theta)\beta}
      \ealigned
    \eequation
    where $\phi_{\min}\in\R{}$ is a lower bound of $\phi(\cdot,\tau_{k_{\tau,\xi}})$ over $\Xcal$ $($by Assumption~\ref{ass.main}$)$.
    \item[(ii)] If $\{\beta_k\}_{k\geq k_{\tau,\xi}+1}$ satisfies $\sum_{k=k_{\tau,\xi}+1}^\infty \beta_k = \infty$ and $\sum_{k=k_{\tau,\xi}+1}^\infty \beta_k^2 < \infty$, then
    \bequation\label{eq:main_thm_diminishing_stepsize}
      \lim_{K\to\infty} \E_{\tau,\low}\left[\tfrac{1}{\sum_{j=k_{\tau,\xi}+1}^{k_{\tau,\xi}+K}\beta_j} \sum_{j=k_{\tau,\xi}+1}^{k_{\tau,\xi}+K} \beta_j\Delta l(x_j,\tau_{k_{\tau,\xi}},\nabla f_j,\dTrue_j)\right] = 0.
    \eequation
  \eitemize
\etheorem
\bproof
  By the definition of $\Abar$, the fact that $\{\beta_k\}\subset (0,1]$, and line~\ref{line:stochastic_choose_stepsize} of Algorithm~\ref{alg:stochastic_sqp_adaptive}, it follows that $\alpha_k \in [\Abar\beta_k,(\Abar + \theta)\beta_k]$ for all $k\geq k_{\tau,\xi} + 1$.  It follows from this fact, $\Delta l(x_k,\tau_{k_{\tau,\xi}},\nabla f_k,\dTrue_k) > 0$ (see~\eqref{Deltal-true}), Lemmas~\ref{lm:stochastic_merit_decrease_upper_bound}, \ref{lm:stochastic_gd_diff}, \ref{lm:Thetak_upper_bound}, \ref{lem:Jd_Jv_diff}, and \ref{lm:stochastic_Delta_q}, and the fact that $\{\beta_k\}\subset (0,1]$ that, for all $k\geq k_{\tau,\xi} + 1$, one finds
  \bequation\label{main:both}
    \baligned
      &\ \E_{k,\tau,\low}[\phi(x_k + \alpha_k d_k,\tau_{k_{\tau,\xi}})] - \phi(x_k,\tau_{k_{\tau,\xi}}) \\
  \leq&\ \E_{k,\tau,\low}[ -\alpha_k\Delta l(x_k,\tau_{k_{\tau,\xi}},\nabla f_k,\dTrue_k) + \alpha_k\tau_{k_{\tau,\xi}}\nabla f_k^T(d_k - \dTrue_k) ] \\
      &\ + (1-\eta)\E_{k,\tau,\low}[ \alpha_k\beta_k\Delta l(x_k,\tau_{k_{\tau,\xi}},g_k,d_k)] \\
      &\ + \E_{k,\tau,\low}[\alpha_k(\|c_k + J_kd_k\| - \|c_k + J_kv_k\|)] \\
  \leq&\ -\Abar\beta_k\Delta l(x_k,\tau_{k_{\tau,\xi}},\nabla f_k,\dTrue_k) + (\omega_8 + \omega_9)\beta_k^2 \\
      &\ + (1-\eta)(\Abar+\theta)\beta_k^2\E_{k,\tau,\low}[ \Delta l(x_k,\tau_{k_{\tau,\xi}},g_k,d_k)] \\
  \leq&\ (- \Abar\beta_k + (1-\eta)(\Abar + \theta)\beta_k^2)\Delta l(x_k,\tau_{k_{\tau,\xi}},\nabla f_k,\dTrue_k) + \beta_k^2\Mbar \\
     =&\ -\beta_k\big(\Abar - (1-\eta)(\Abar + \theta)\beta_k\big)\Delta l(x_k,\tau_{k_{\tau,\xi}},\nabla f_k,\dTrue_k) + \beta_k^2\Mbar.
    \ealigned
  \eequation
  Let us now consider the two cases in the theorem one at a time.

  \textbf{Case (i).} By the definition of $\beta$, it follows by taking total expectation of~\eqref{main:both} (namely, expectation defined in~\eqref{def:Etaulow}) that for each $k\geq k_{\tau,\xi} + 1$ one has
  \bequationNN
    \baligned
      &\ \E_{\tau,\low}[\phi(x_k + \alpha_k d_k,\tau_{k_{\tau,\xi}})] - \E_{\tau,\low}[\phi(x_k,\tau_{k_{\tau,\xi}})] \\
 \leq &\ -\beta(\Abar - (1-\eta)(\Abar + \theta)\beta)\E_{\tau,\low}[\Delta l(x_k,\tau_{k_{\tau,\xi}},\nabla f_k,\dTrue_k)] + \beta^2\Mbar.
    \ealigned
  \eequationNN
  Summing this inequality over $j\in\{k_{\tau,\xi}+1,\ldots,k_{\tau,\xi}+K\}$ shows that
  \bequationNN
    \baligned
      &\ \phi_{\min} - \E_{\tau,\low}[\phi(x_{k_{\tau,\xi}+1},\tau_{k_{\tau,\xi}})] \\
  \leq&\ \E_{\tau,\low}[\phi(x_{k_{\tau,\xi}+K+1},\tau_{k_{\tau,\xi}})] - \E_{\tau,\low}[\phi(x_{k_{\tau,\xi}+1},\tau_{k_{\tau,\xi}})] \\
  \leq&\ -\beta (\Abar - (1-\eta)(\Abar + \theta)\beta)\E_{\tau,\low} \left[ \sum_{j=k_{\tau,\xi}+1}^{k_{\tau,\xi}+K} \Delta l(x_j,\tau_{k_{\tau,\xi}},\nabla f_j,\dTrue_j) \right] + K\beta^2\Mbar,
    \ealigned
  \eequationNN
  which after rearrangement shows that~\eqref{eq:main_thm_fixed_stepsize} holds, as desired.

  \textbf{Case (ii).}  Given the definition of $\{\beta_k\}$, let us assume without loss of generality that $\beta_k \leq \Abar/\big(2(1-\eta)(\Abar + \theta)\big)$ for all $k\geq k_{\tau,\xi}+1$, which implies that $\Abar - (1-\eta)(\Abar+\theta)\beta_k \geq \tfrac{1}{2}\Abar$ for all $k \geq k_{\tau,\xi}+1$.  Using this fact, taking total expectation of~\eqref{main:both} (namely, expectation defined in~\eqref{def:Etaulow}), and using~\eqref{Deltal-true} it holds that
  \bequationNN
    \baligned
      &\ \E_{\tau,\low}[\phi(x_k + \alpha_k d_k,\tau_{k_{\tau,\xi}})] - \E_{\tau,\low}[\phi(x_k,\tau_{k_{\tau,\xi}})] \\
  \leq&\ -\tfrac{1}{2}\beta_k\Abar\E_{\tau,\low}[\Delta l(x_k,\tau_{k_{\tau,\xi}},\nabla f_k,\dTrue_k)] + \beta_k^2\Mbar.
    \ealigned
  \eequationNN
  Summing this inequality over $j\in\{k_{\tau,\xi}+1,\ldots, k_{\tau,\xi}+K\}$ shows that
  \bequationNN
    \baligned
      &\ \phi_{\min} - \E_{\tau,\low}[\phi(x_{k_{\tau,\xi}+1},\tau_{k_{\tau,\xi}})] \\
  \leq&\ \E_{\tau,\low}[\phi(x_{k_{\tau,\xi}+K+1},\tau_{k_{\tau,\xi}})] - \E_{\tau,\low}[\phi(x_{k_{\tau,\xi}+1},\tau_{k_{\tau,\xi}})] \\
  \leq&\ -\tfrac{1}{2}\Abar\E_{\tau,\low} \left[\sum_{j=k_{\tau,\xi}+1}^{k_{\tau,\xi}+K} \beta_j\Delta l(x_j,\tau_{k_{\tau,\xi}},\nabla f_j,\dTrue_j) \right] + \Mbar\sum_{j=k_{\tau\xi}+1}^{k_{\tau,\xi}+K}\beta_j^2,
    \ealigned
  \eequationNN
  which after rearrangement and taking limits proves that \eqref{eq:main_thm_diminishing_stepsize} holds.
\eproof

\section{Numerical Results}\label{sec.numerical}

In this section, we demonstrate the performance of a Matlab implementation of Algorithm~\ref{alg:stochastic_sqp_adaptive} for solving (i) a subset of the CUTEst collection of test problems~\cite{GoulOrbaToin15} and (ii) two optimal control problems from~\cite{HintItoKuni03}.  The goal of our testing is to demonstrate the computational benefits of using inexact subproblem solutions obtained based on our termination tests from Section~\ref{sec:direction}. 

\subsection{Iterative solvers}\label{sec:solvers}

To obtain the normal direction $v_k$ as an inexact solution of~\eqref{prob.normal}, we applied the conjugate gradient (CG) method to $J_k^TJ_kv = -J_k^Tc_k$.  Denoting the $t$th CG iterate as $v_{k,t}$, where $v_{k,0} = 0$, the method sets $v_k \gets v_{k,t}$, where~$t$ is the first CG iteration such that $\|J_k^TJ_kv_{k,t} + J_k^Tc_k\| \leq \max\{0.1\|J_k^Tc_k\|,10^{-10}\}$.  The properties of the CG method as a Krylov subspace method ensure that $v_{k,t}\in\Range(J_k^T)$ for all $t\in\N{}$ (in exact arithmetic); hence, $v_k\in\Range(J_k^T)$.

To obtain the tangential direction $u_k$ and associated dual search direction $\delta_k$, we applied the minimum residual (MINRES) method, namely, the implementation from~\cite{ChoiPaigSaun11,PaigSaun75}, to the linear system
\bequation\label{eq:stochastic_system_numerics}
  \bbmatrix H_k & J_k^T \\ J_k & 0 \ebmatrix \bbmatrix u \\ \delta \ebmatrix = - \bbmatrix g_k + H_kv_k + J_k^Ty_k \\ 0 \ebmatrix.
\eequation
(We discuss our choice of $H_k$ along with each set of experiments.)  Letting $(u_{k,t},\delta_{k,t})$ denote the $t$th MINRES iterate, where $(u_{k,0},\delta_{k,0}) = (0,0)$, the method sets $(u_k,\delta_k) \gets (u_{k,t},\delta_{k,t})$ where~$t$ is the first MINRES iteration such that, for some $\kappa\in(0,1)$, \bequation\label{eq:additional_condition}
  \left\|\bbmatrix \rho_{k,t} \\ r_{k,t} \ebmatrix \right\|_{\infty} \leq \max \left\{ \kappa \left\|\bbmatrix g_k + H_kv_k + J_k^Ty_k \\ 0 \ebmatrix \right\|_{\infty}, 10^{-12} \right\}
\eequation
and Termination Test 1 and/or 2 holds.  (Recall the definition of $(\rho_{k,t},r_{k,t})$ in~\eqref{eq:stochastic_system_iterative}.)  The choice of $\kappa\in(0,1)$ is discussed along with each set of experiments.

\subsection{Choosing the step size}

Algorithm~\ref{alg:stochastic_sqp_adaptive} (see line~\ref{line:stochastic_choose_stepsize}) stipulates that the step size $\alpha_k$ chosen for the $k$th iteration satisfies $\alpha_k\in[\alphamin_k,\alphamax_k]$.  Keeping in mind that the inequalities $\alphamin_k \leq \alphasuff_k \leq 1$ and $\alphasuff_k \leq \alphavarphi_k$ (see Lemma~\ref{lm:stepsize_well_defined}) hold, we take advantage of this flexibility in choosing the step size by defining
\bequationNN
  \alpha_k \gets \bcases \min\{\alphasuff_k,\alphamin_k+\theta\beta_k^2\} & \text{if $\alphasuff_k = 1$} \\
  \alphamin_k+\theta\beta_k^2 & \text{if $\alphamin_k+\theta\beta_k^2 \leq \alphasuff_k < 1$} \\
  (1.1)^{t_k}\alphasuff_k & \text{if $\alphasuff_k < \min\{\alphamin_k+\theta\beta_k^2, 1\}$,} \ecases
\eequationNN 
where $t_k$ is the largest value of $t \in \N{}$ such that
\bequationNN
  (1.1)^t\alphasuff_k \leq \min\{\alphavarphi_k,\alphamin_k+\theta\beta_k^2\} \equiv \alphamax_k \ \text{and} \ (1.1)^{t-1}\alphasuff_k < 1.
\eequationNN
When $\alphasuff_k < 1$, this strategy allows for the possibility that step sizes larger than $\min\{\alphasuff_k , \alphamin_k + \theta\beta_k^2\}$ be taken (namely, they can be as large as $\min\{\alphavarphi_k , \alphamin_k + \theta\beta_k^2\}$).  This led to better performance while still having a rule that satisfies the requirements of our analysis.  We do not explicitly compute $\alphavarphi_k$ in our code.  Instead, we can verify directly whether $(1.1)^t\alphasuff_k \leq \alphavarphi_k$ (as needed above) since this is ensured by checking whether $\varphi((1.1)^t\alphasuff_k) \leq 0$, which is easily checked by the code.

\subsection{Algorithm variants tested}

To test the utility of using inexact subproblem solutions in Algorithm~\ref{alg:stochastic_sqp_adaptive}, we consider two algorithm variants that we refer to as \SISQO{} and \SISQOexact{}.  The variant \SISQO{} is Algorithm~\ref{alg:stochastic_sqp_adaptive} with inexact solutions computed as described in Section~\ref{sec:solvers} with a relatively \emph{large} value for $\kappa$ in~\eqref{eq:additional_condition}.  On the other hand, the variant \SISQOexact{} is identical to \SISQO{} with the exception that it uses a relatively \emph{small} value for $\kappa$ in~\eqref{eq:additional_condition}.  We specify the values of $\kappa \in (0,1)$ used along with each of our tests in Sections~\ref{sec.CUTE} and \ref{sec.control}.

Our reason for comparing these two variants is to focus attention on the numerical gains obtained as a result of using inexact subproblem solutions. For this reason, we allow both variants to use the same computation for the normal step, thus allowing any numerical gains to be directly attributed to the inexact tangential step computation.  Although other variants could be tested (e.g., allowing the normal step computation to   differ as well) we prefer the approach described above since it limits the variation attributable to the different calculations in the SISQO framework.

\subsection{Metrics used for comparison}\label{sec:metrics}

Our metrics of interest are feasibility and stationarity. Specifically, for any run of \SISQO{}, we terminate with $\xinexact \gets x_k$, where $k \in \N{}$ is the first iteration such that $\|c(x_k)\|_{\infty} \leq 10^{-6}$ and $\|\nabla f_k + J_k^T y_{k,{\rm ls}}\|_{\infty} \leq 10^{-2}$, where $y_{k,{\rm ls}}$ is the least-square multiplier at $x_k$.  (The computations of $\nabla f_k$ and $y_{k,{\rm ls}}$ are not required by our algorithm in general; they were computed in our experiments merely for the purpose of being able to determine an accurate measure of stationarity at $x_k$.)  This allows us to associate with each run of \SISQO{} the two measures
\bequationNN
\baligned
\text{error}_\text{feasibility}(\SISQO{}) &= \|c(\xinexact)\|_{\infty}
\ \ \text{and} \ \  \\
\text{error}_\text{stationarity}(\SISQO{}) &= \|\nabla f(\xinexact) + J(\xinexact)^T \yinexact\|_{\infty},
\ealigned
\eequationNN
where $\yinexact\in\R{m}$ is the least-square multiplier at $\xinexact$.  We use the total number of MINRES iterations performed by \SISQO{} as a budget for the number of MINRES iterations performed by \SISQOexact{}; no other termination condition is used for \SISQOexact{}.  Upon termination of \SISQOexact{}, we define $\xexact$ in the following manner: If an iterate is computed with $\|c(x_k)\|_{\infty} \leq 10^{-6}$, then $\xexact$ is chosen as the iterate with smallest stationarity measure among those satisfying this tolerance for the feasibility measure; otherwise, $\xexact$ is chosen as the iterate with the smallest feasibility measure. In any case, once $\xexact$ is determined, we proceed to compute the least-square multiplier $\yexact$ at $\xexact$, then define 
\bequationNN
\baligned
\text{error}_\text{feasibility}(\SISQOexact{}) &= \|c(\xexact)\|_{\infty}
\ \ \text{and} \ \  \\
\text{error}_{\text{stationarity}}(\SISQOexact{}) &= \|\nabla f(\xexact) + J(\xexact)^T \yexact\|_{\infty}.
\ealigned
\eequationNN
These are the metrics that we use in the next two subsections.

\subsection{Results on the CUTEst problems}\label{sec.CUTE}

In the CUTEst collection \cite{GoulOrbaToin15}, there are a total of 138 equality constrained problems with $m\leq n$. From these problems, we selected those such that (i) $(n+m)\in [500,10000]$, (ii) the objective function is not constant, (iii) the objective function remained above $-10^{50}$ over the sequences of iterates generated by runs of our algorithm, and (iv) the LICQ was satisfied at all iterates encountered in each run of our algorithm.  This process of elimination resulted in the following $11$ test problems: \texttt{ELEC}, \texttt{LCH}, \texttt{LUKVLE1}, \texttt{LUKVLE3}, \texttt{LUKVLE4}, \texttt{LUKVLE6}, \texttt{LUKVLE7}, \texttt{LUKVLE9}, \texttt{LUKVLE10}, \texttt{LUKVLE13}, and \texttt{ORTHREGC}.

The problems from the CUTEst collection are deterministic, and for the purpose of these experiments we exploited this fact to compute values as needed by our algorithm, including using function evaluations to estimate Lipschitz constants and using a (modified) Hessian of the Lagrangian in each search direction computation, as explained below.  However, we introduced noise into the computation of the objective function gradients.  In particular, we generated stochastic gradients as  $g_k = \Ncal(\nabla f_k,\tfrac{\epsilon_N^2}{n}I)$, where for testing purposes we considered the three noise levels $\epsilon_N \in \{10^{-4},10^{-2},10^{-1}\}$. This particular choice for defining the stochastic gradients ensured that an appropriate value for $M_g$ as indicated in Assumption~\ref{ass.g} would be given by $M_g = \{10^{-8},10^{-4},10^{-2}\}$, corresponding to the values for $\epsilon_N$.

In terms of algorithm parameters, we set $\kappa = 0.1$ for \SISQO{} and $\kappa = 10^{-7}$ for 
\SISQOexact{}. All of the remaining parameters were set identically for the two variants with the following values:
$\tau_{-1} = \sigma_c = \eta = \kappa_v = \kappa_u = 0.1$, $\xi_{-1} = \epsilon_c = 1$, $\epsilon_{\tau} = \epsilon_{\xi} = 0.01$, $\kappa_{\rho} = \kappa_r = 100$, $\epsilon_r = 1-10^{-4}$, $\zeta = 10^{-8}$, $\epsilon_u = 5\times 10^{-9}$, $\sigma_u = 1-10^{-12}$, $\theta = 10^4$, and $\beta_k = 1$ for all $k\in\N{}$. During each iteration $k\in\N{}$, we randomly generated a sample point near $x_k$, then estimated $L_k$ and $\Gamma_k$ using finite differences of the objective gradients and constraint Jacobians between $x_k$ and the sampled point.  These values were used in place of $L$ and $\Gamma$, respectively, in our step size selection.

For this collection of problems, we employed an iterative Hessian modification strategy as proposed in~\cite{ByrdCurtNoce10,CurtNoceWach09}.  Specifically, for all $k \in \N{}$, the matrix $H_k$ is initialized to the true Hessian of the Lagrangian, but may be set ultimately as
\bequationNN
  H_k \gets \iota_k \nabla^2_{xx}\big(f(x)+c(x)^T y\big)|_{(x,y)=(x_k,y_k)} + (1 - \iota_k) I
\eequationNN
with $\iota_k = 10^{-j_k}$, where $j_k$ is the smallest element in $\{0,\dots,10\}$ such that a modification is not triggered.  If a modification is triggered at $\iota_k = 10^{-10}$, then the algorithm sets $H_k \gets I$ to guarantee that no further modifications are required.

\begin{figure}[ht]
  \centering
  \includegraphics[width=0.42\textwidth,clip=true,trim=20 5 90 50]{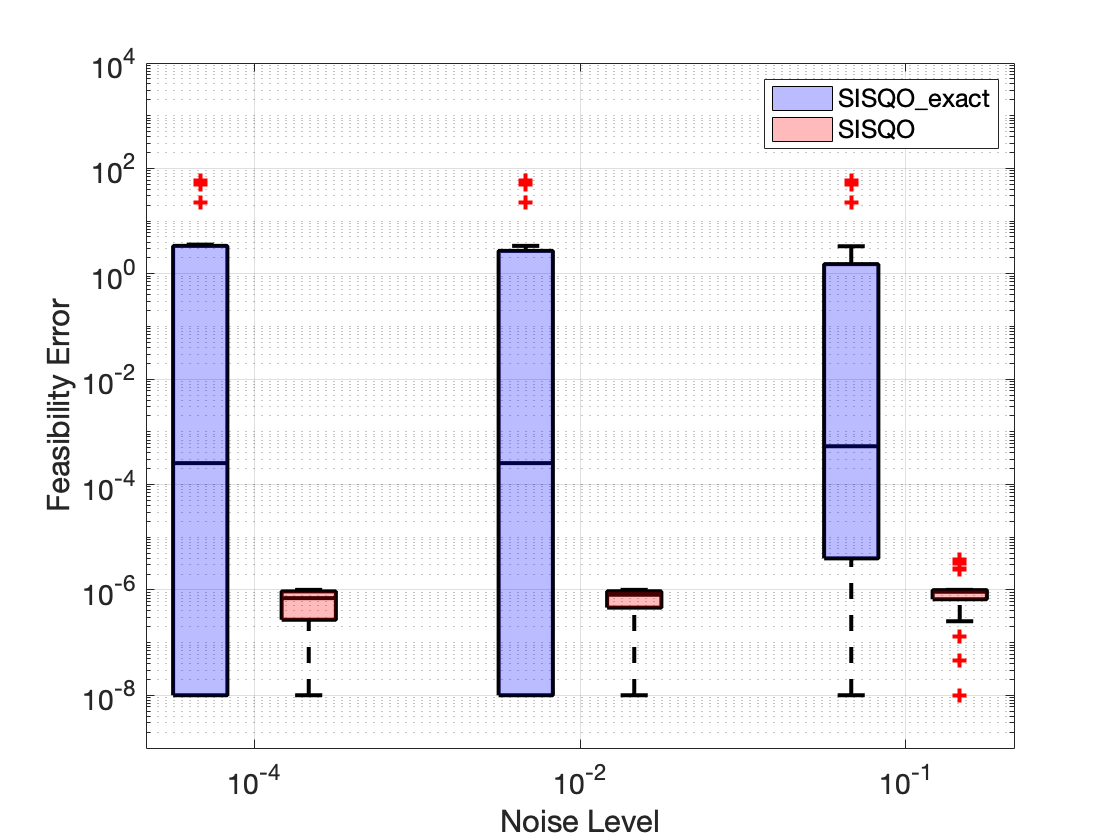}\quad 
 \includegraphics[width=0.42\textwidth,clip=true,trim=20 5 90 50]{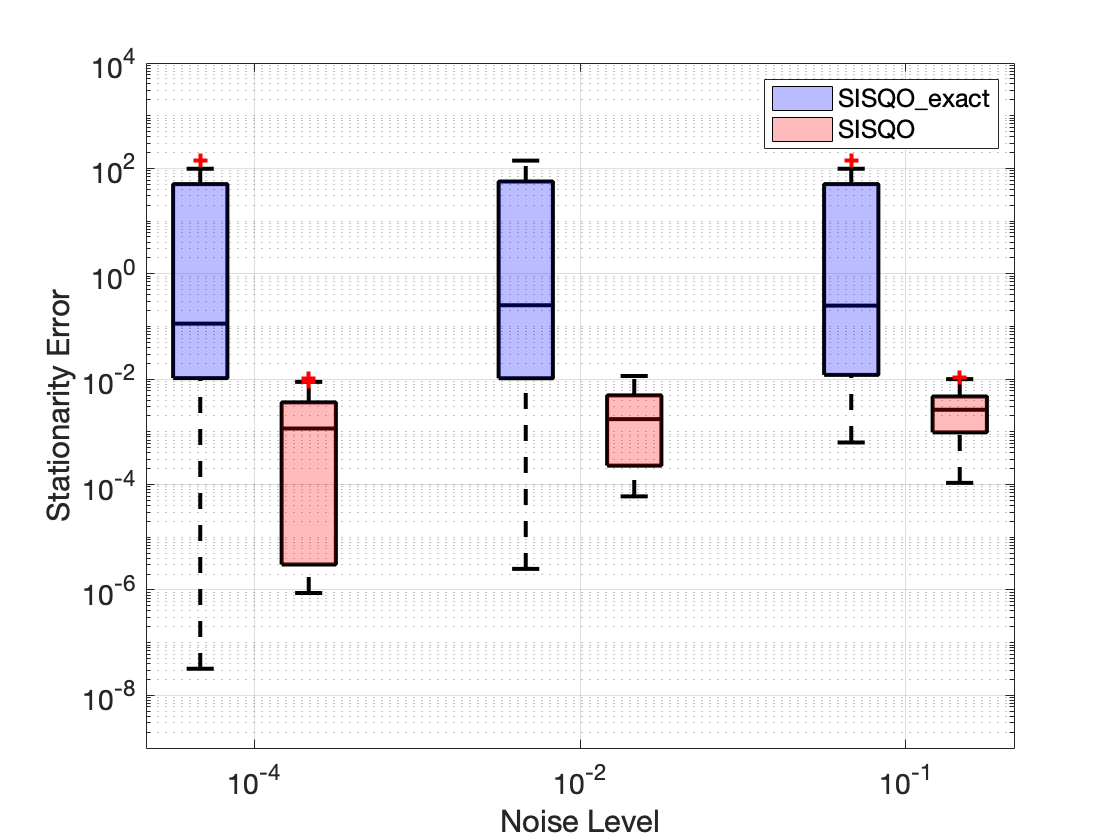}
  \caption{Box plots on CUTEst problems for feasibility (left) and stationarity (right).}
  \label{fig.perf_cutest}
\end{figure}

For each test problem, we ran both \SISQO{} and \SISQOexact{} five times, and for each computed the resulting feasibility and stationarity errors as described in Section~\ref{sec:metrics}. The results are shown in the form of box plots in Figure~\ref{fig.perf_cutest}.

From Figure~\ref{fig.perf_cutest}, one finds that \SISQO{} performs better than \SISQOexact{} in terms of both feasibility and stationarity errors.  Also, in general, \SISQO{} achieves smaller feasibility and stationarity errors for smaller noise levels, which may be expected due to the fact that these experiments are run with constant $\{\beta_k\}$.

\subsection{Results on optimal control problems}\label{sec.control}

In our second set of experiments, we considered two optimal control problems motivated by those in~\cite{HintItoKuni03}.  In particular, we modified the problems to have equality constraints only and finite sum objective functions.  Specifically, given a domain $\Xi\in\R{2}$, a constant $N\in\N{}_{>0}$, reference functions $\wbar_{ij}\in L^2(\Xi)$ and 
$\zbar_{ij}\in L^2(\Xi)$ 
for $(i,j)\in\{1,\ldots,N\}\times\{1,\ldots,N\}$, and a regularization parameter $\lambda\in\R{}_{>0}$, we first considered the problem
\bequation\label{eq:OC_prob_1}
  \baligned
    \min_{w,z}\ & \frac{1}{N^2}\sum_{i=1}^N\sum_{j=1}^N (\tfrac{1}{2} \|w - \wbar_{ij}\|_{L^2(\Xi)}^2 + \tfrac{\lambda}{2} \|z - \zbar_{ij}\|_{L^2(\Xi)}^2) \\
    \text{s.t. }& -\Delta w = z \ \text{in} \ \Xi, \ \text{and} \ w = 0 \ \text{on}\ \partial\Xi.
  \ealigned
\eequation
Second, with the same notation but $\zbar_{ij}\in L^2(\partial \Xi)$, we also considered
\bequation\label{eq:OC_prob_2}
  \baligned
    \min_{w,z}\ & \frac{1}{N^2}\sum_{i=1}^N\sum_{j=1}^N (\tfrac{1}{2} \|w - \wbar_{ij}\|_{L^2(\Xi)}^2 + \tfrac{\lambda}{2} \|z - \zbar_{ij}\|_{L^2(\partial\Xi)}^2) \\
    \text{s.t. }& -\Delta w + w = 0 \ \text{in} \ \Xi, \ \text{and} \ \tfrac{\partial w}{\partial p} = z \ \text{on} \ \partial\Xi,
  \ealigned
\eequation
where $p$ represents the unit outer normal to $\Xi$ along $\partial \Xi$.
As reference functions for both problems, we chose for all $(i,j)\in\{1,\ldots,N\}\times\{1,\ldots,N\}$ the following:
\bequation\label{eq:reference_functions}
  \zbar_{ij} = 0\ \text{and}\ \wbar_{ij}(x_1,x_2) = \sin((4 + \tfrac{\epsilon_N}{\epsilon_S}  (i-\tfrac{N+1}{2})) x_1) + \cos((3 + \tfrac{\epsilon_N}{\epsilon_S} (j-\tfrac{N+1}{2})) x_2)
\eequation
for some $(\epsilon_S,\epsilon_N)\in\R{}_{>0}\times\R{}_{>0}$.  We selected the following values for the above constants: 
$N=3$, $\lambda = 10^{-5}$, $\epsilon_S = \sqrt{15}$, and $\epsilon_N \in \{10^{-4},10^{-2},10^{-1}\}$.  Since the objective functions of \eqref{eq:OC_prob_1} and \eqref{eq:OC_prob_2} are finite sums, to generate stochastic gradients as unbiased estimates of the true gradient, we first uniformly generated random $(i,j)\in \{1,\ldots,N\} \times \{1,\ldots,N\}$, then computed the gradient corresponding to the $(i,j)$th term in the objective function.  We note that with the above choice of parameters, it follows that an appropriate value for $M_g$ in Assumption~\ref{ass.g} is given by $M_g \approx \{10^{-8},10^{-4},10^{-2}\}$ to correspond, respectively, to the above values for $\epsilon_N$.

Since the optimal control problems have a quadratic objective function and linear constraints, we used the exact second derivative matrix $H_k = \diag(I,\lambda I)$ for all $k\in\N{}$. For this choice, the curvature condition on $H_k$ in Assumption~\ref{ass.H} is trivially satisfied.

In terms of algorithm parameters, we set $\kappa = 10^{-4}$ for \SISQO{} and $\kappa = 10^{-7}$ for \SISQOexact{}. All of the remaining parameters were set identically for the two variants in the same manner as in the previous section with the following exceptions: $\tau_{-1} = 10^{-4}$, $\eta = 0.5$, and $L_k = 1$ and $\Gamma_k = 0$ for all $k\in\N{}$, which are valid choices since the objective functions are quadratic and the constraints are linear. 

For each of the two optimal control problems in \eqref{eq:OC_prob_1} and \eqref{eq:OC_prob_2}, we ran both \SISQO{} and \SISQOexact{} ten times, then computed their average feasibility and stationarity errors as described in Section~\ref{sec:metrics}.  In Table~\ref{table:OC_prob_1} and Table~\ref{table:OC_prob_2}, we report
these average values as well as the average number of iterations performed by  Algorithm~\ref{alg:stochastic_sqp_adaptive} before termination (``iterations") and number of MINRES iterations (``MINRES iterations"), with the latter discussed in Section~\ref{sec:solvers}. The results are given in Table~\ref{table:OC_prob_1} and Table~\ref{table:OC_prob_2} for problem~\eqref{eq:OC_prob_1} and problem~\eqref{eq:OC_prob_2}, respectively.  One can observe that despite performing more ``outer'' iterations on average, \SISQO{} outperforms \SISQOexact{} due to the fact that it requires fewer overall linear system solver iterations on average in order to attain better average feasibility and stationarity errors.

\begin{table}[ht]
\centering
\small
\begin{tabular}{ |l|c|c|c|c|c| } 
 \hline
 strategy & $\epsilon_N$ & 
 \begin{tabular}{c} feasibility \\[-0.2em] error \end{tabular} & 
 \begin{tabular}{c} stationarity \\[-0.2em] error \end{tabular} &
 \begin{tabular}{c} MINRES \\[-0.2em] iterations \end{tabular} &
 iterations \\ 
 \hline 
 \hline 
 \SISQO{} & $10^{-4}$ & $2.41\times 10^{-7}$ & $1.76\times 10^{-5}$ & 55117 & 8.9 \\ 
 \hline
 \SISQOexact{} & $10^{-4}$ & $3.86\times 10^{-5}$ & $4.05\times 10^{-5}$ & 55117 & 6.9 \\
 \hline\hline
 \SISQO{} & $10^{-2}$ & $4.14\times 10^{-7}$ & $2.09\times 10^{-3}$ & 60894 & 8.8 \\ 
 \hline
 \SISQOexact{} & $10^{-2}$ & $3.46\times 10^{-5}$ & $1.95\times 10^{-3}$ & 60894 & 6.8 \\ 
 \hline\hline
 \SISQO{} & $10^{-1}$ & $3.43\times 10^{-7}$ & $5.15\times 10^{-3}$ & 93634 & 12.3 \\ 
 \hline
 \SISQOexact{} & $10^{-1}$ & $2.36\times 10^{-6}$ & $1.68\times10^{-2}$ & 93634 & 10 \\ 
 \hline
\end{tabular}
\caption{Numerical results for problem~\eqref{eq:OC_prob_1} averaged over ten independent runs.}
\label{table:OC_prob_1}
\end{table}

\begin{table}[ht]
\centering
\small
\begin{tabular}{ |l|c|c|c|c|c| } 
 \hline
 strategy & $\epsilon_N$ &
 \begin{tabular}{c} feasibility \\[-0.2em] error \end{tabular} & 
 \begin{tabular}{c} stationarity \\[-0.2em] error \end{tabular} &
 \begin{tabular}{c} MINRES \\[-0.2em] iterations \end{tabular} &
iterations \\
 \hline 
 \hline
 \SISQO{} & $10^{-4}$ & $3.29\times 10^{-7}$ & $2.35\times 10^{-5}$ & 91478 & 9.9 \\ 
 \hline
 \SISQOexact{} & $10^{-4}$ & $5.44\times 10^{-4}$ & $5.46\times 10^{-4}$ & 91478 & 7.1 \\
 \hline\hline
 \SISQO{} & $10^{-2}$ & $2.90\times 10^{-7}$ & $2.07\times 10^{-3}$ & 99921 & 10 \\ 
 \hline
 \SISQOexact{} & $10^{-2}$ & $5.71\times 10^{-5}$ & $2.37\times 10^{-3}$ & 99921 & 7.6 \\ 
 \hline\hline
 \SISQO{} & $10^{-1}$ & $1.68\times 10^{-7}$ & $3.88\times 10^{-4}$ & 158825 & 14.5 \\ 
 \hline
 \SISQOexact{} & $10^{-1}$ & $1.31\times 10^{-5}$ & $2.58\times 10^{-2}$ & 158825 & 11.1 \\
 \hline
\end{tabular}
\caption{Numerical results for  problem \eqref{eq:OC_prob_2} averaged over ten independent runs.}
\label{table:OC_prob_2}
\end{table}

\section{Conclusion}\label{sec.conclusion}

We have proposed, analyzed, and tested an \emph{inexact} stochastic SQP algorithm for solving stochastic optimization problems involving deterministic, smooth, nonlinear equality constraints.  We proved a convergence guarantee (in expectation) for our algorithm that is comparable to that proved for the \emph{exact} stochastic SQP method recently presented in \cite{BeraCurtRobiZhou21}, which in turn is comparable to that known for the stochastic gradient in unconstrained settings \cite{BottCurtNoce18}.  Our Matlab implementation, \SISQO{}, illustrated the benefits of allowing inexact step computation for solving problems from the CUTEst collection~\cite{GoulOrbaToin15} as well as two optimal control problems.

\bibliographystyle{plain}
\bibliography{references}

\begin{thebibliography}{10}

\bibitem{AchiHeldTamaAbbe17}
Joshua Achiam, David Held, Aviv Tamar, and Pieter Abbeel.
\newblock {Constrained policy optimization}.
\newblock In {\em Proceedings of the 34th International Conference on Machine
  Learning-Volume 70}, pages 22--31, 2017.

\bibitem{BeraCurtONeiRobi21}
Albert~S. Berahas, Frank~E. Curtis, Michael~J. O'Neill, and Daniel~P. Robinson.
\newblock A stochastic sequential quadratic optimization algorithm for
  nonlinear equality constrained optimization with rank-deficient {J}acobians.
\newblock {\em arXiv preprint arXiv:2106.13015}, 2021.

\bibitem{BeraCurtRobiZhou21}
Albert~S. Berahas, Frank~E. Curtis, Daniel~P. Robinson, and Baoyu Zhou.
\newblock Sequential quadratic optimization for nonlinear equality constrained
  stochastic optimization.
\newblock {\em SIAM Journal on Optimization}, 31(2):1352--1379, 2021.

\bibitem{BiroGhat03}
G.~Biros and O.~Ghattas.
\newblock {Inexactness Issues in the Lagrange-Newton-Krylov-Schur Method for
  PDE-constrained Optimization}.
\newblock In L.~T. Biegler, O.~Ghattas, M.~Heinkenschloss, and B.~{Van Bloemen
  Waanders}, editors, {\em Large-Scale PDE-Constrained Optimization}, pages
  93--114, New York, NY, USA, 2003. Springer.

\bibitem{BottCurtNoce18}
L{\'e}on Bottou, Frank~E Curtis, and Jorge Nocedal.
\newblock Optimization methods for large-scale machine learning.
\newblock {\em Siam Review}, 60(2):223--311, 2018.

\bibitem{ByrdCurtNoce08}
Richard~H. Byrd, Frank~E. Curtis, and Jorge Nocedal.
\newblock An inexact {SQP} method for equality constrained optimization.
\newblock {\em {SIAM Journal on Optimization}}, 19(1):351--369, 2008.

\bibitem{ByrdCurtNoce10}
Richard~H Byrd, Frank~E Curtis, and Jorge Nocedal.
\newblock An inexact {Newton} method for nonconvex equality constrained
  optimization.
\newblock {\em Mathematical programming}, 122(2):273, 2010.

\bibitem{ChatChenMaasCarr16}
Nilanjan Chatterjee, Yi-Hau Chen, Paige Maas, and Raymond~J. Carroll.
\newblock Constrained maximum likelihood estimation for model calibration using
  summary-level information from external big data sources.
\newblock {\em Journal of the American Statistical Association},
  111(513):107--117, 2016.

\bibitem{ChenTungVeduMori18}
Changan Chen, Frederick Tung, Naveen Vedula, and Greg Mori.
\newblock {Constraint-aware deep neural network compression}.
\newblock In {\em Proceedings of the European Conference on Computer Vision
  (ECCV)}, pages 400--415, 2018.

\bibitem{ChoiPaigSaun11}
Sou-Cheng~T Choi, Christopher~C Paige, and Michael~A Saunders.
\newblock {MINRES-QLP: A Krylov} subspace method for indefinite or singular
  symmetric systems.
\newblock {\em SIAM Journal on Scientific Computing}, 33(4):1810--1836, 2011.

\bibitem{Cour43}
R~Courant.
\newblock Variational methods for the solution of problems of equilibrium and
  vibrations.
\newblock {\em Bulletin of the American Mathematical Society}, 49(1):1--23,
  1943.

\bibitem{CurtNoceWach09}
Frank~E Curtis, Jorge Nocedal, and Andreas W{\"a}chter.
\newblock A matrix-free algorithm for equality constrained optimization
  problems with rank-deficient {J}acobians.
\newblock {\em SIAM Journal on Optimization}, 20(3):1224--1249, 2009.

\bibitem{Flet13}
Roger Fletcher.
\newblock {\em {Practical Methods of Optimization}}.
\newblock John Wiley \& Sons, 2013.

\bibitem{Geye91}
Charles~J. Geyer.
\newblock Constrained maximum likelihood exemplified by isotonic convex
  logistic regression.
\newblock {\em Journal of the American Statistical Association},
  86(415):717--724, 1991.

\bibitem{GoulOrbaToin15}
Nicholas I.~M. Gould, Dominique Orban, and Philippe~L. Toint.
\newblock {CUTEst}: a constrained and unconstrained testing environment with
  safe threads for mathematical optimization.
\newblock {\em Computational Optimization and Applications}, 60:545--557, 2015.

\bibitem{HanMang79}
S-P Han and Olvi~L Mangasarian.
\newblock Exact penalty functions in nonlinear programming.
\newblock {\em Mathematical programming}, 17(1):251--269, 1979.

\bibitem{Han77}
Shih-Ping Han.
\newblock A globally convergent method for nonlinear programming.
\newblock {\em Journal of optimization theory and applications},
  22(3):297--309, 1977.

\bibitem{HeinRidz08b}
M.~Heinkenschloss and D.~Ridzal.
\newblock {An Inexact Trust-Region SQP Method with Applications to
  PDE-Constrained Optimization}.
\newblock In K.~Kunisch, G.~Of, and O.~Steinbach, editors, {\em Numerical
  Mathematics and Advanced Applications: Proceedings of ENUMATH 2007, the 7th
  European Conference on Numerical Mathematics and Advanced Applications, Graz,
  Austria}, pages 613--620. Springer, 2008.

\bibitem{HeinVice02}
M.~Heinkenschloss and L.~N. Vicente.
\newblock {Analysis of Inexact Trust-Region SQP Algorithms}.
\newblock {\em SIAM Journal on Optimization}, 12(2):283--302, 2002.

\bibitem{HintItoKuni03}
Michael Hinterm{\"u}ller, Kazufumi Ito, and Karl Kunisch.
\newblock The primal-dual active set strategy as a semismooth {Newton} method.
\newblock {\em SIAM Journal on Optimization}, 13(3):865--888, 2003.

\bibitem{KaplTich98}
Alexander Kaplan and Rainer Tichatschke.
\newblock Proximal point methods and nonconvex optimization.
\newblock {\em Journal of global Optimization}, 13(4):389--406, 1998.

\bibitem{KumaSoumMhamHara18}
Soumava Kumar~Roy, Zakaria Mhammedi, and Mehrtash Harandi.
\newblock {Geometry aware constrained optimization techniques for deep
  learning}.
\newblock In {\em Proceedings of the IEEE Conference on Computer Vision and
  Pattern Recognition}, pages 4460--4469, 2018.

\bibitem{MarqSalzFua17}
Pablo M{\'a}rquez-Neila, Mathieu Salzmann, and Pascal Fua.
\newblock Imposing hard constraints on deep networks: Promises and limitations.
\newblock {\em arXiv preprint arXiv:1706.02025}, 2017.

\bibitem{NaAnitKola21}
Sen Na, Mihai Anitescu, and Mladen Kolar.
\newblock An adaptive stochastic sequential quadratic programming with
  differentiable exact augmented {Lagrangians}.
\newblock {\em arXiv preprint arXiv:2102.05320}, 2021.

\bibitem{NandPathAbhiSing19}
Yatin Nandwani, Abhishek Pathak, and Parag Singla.
\newblock {A primal-dual formulation for deep learning with constraints}.
\newblock In {\em Advances in Neural Information Processing Systems}, pages
  12157--12168, 2019.

\bibitem{NoceWrig06}
J.~Nocedal and S.~J. Wright.
\newblock {\em {Numerical Optimization}}.
\newblock Springer Series in Operations Research. Springer, New York, NY, USA,
  second edition, 2006.

\bibitem{PaigSaun75}
Christopher~C Paige and Michael~A Saunders.
\newblock Solution of sparse indefinite systems of linear equations.
\newblock {\em SIAM journal on numerical analysis}, 12(4):617--629, 1975.

\bibitem{Powe78b}
M.~J.~D. Powell.
\newblock A fast algorithm for nonlinearly constrained optimization
  calculations.
\newblock In {\em {Numerical Analysis}}, Lecture Notes in Mathematics, pages
  144--157. Springer, Berlin, 1978.

\bibitem{Powe78a}
Michael~JD Powell.
\newblock Algorithms for nonlinear constraints that use {Lagrangian} functions.
\newblock {\em Mathematical programming}, 14(1):224--248, 1978.

\bibitem{RaviDinhLokhSing19}
Sathya~N Ravi, Tuan Dinh, Vishnu~Suresh Lokhande, and Vikas Singh.
\newblock {Explicitly imposing constraints in deep networks via conditional
  gradients gives improved generalization and faster convergence}.
\newblock In {\em Proceedings of the AAAI Conference on Artificial
  Intelligence}, volume~33, pages 4772--4779, 2019.

\bibitem{Rock76}
R~Tyrrell Rockafellar.
\newblock Monotone operators and the proximal point algorithm.
\newblock {\em SIAM journal on control and optimization}, 14(5):877--898, 1976.

\bibitem{RuthHabe20}
Lars Ruthotto and Eldad Haber.
\newblock {Deep Neural Networks Motivated by Partial Differential Equations}.
\newblock {\em Journal of Mathematical Imaging and Vision}, 62:352--364, 2020.

\bibitem{ShapDentRusz14}
Alexander Shapiro, Darinka Dentcheva, and Andrzej Ruszczy{\'n}ski.
\newblock {\em Lectures on stochastic programming: modeling and theory}.
\newblock SIAM, 2014.

\bibitem{SherRaguMoreAdamThan19}
Sheroze Sheriffdeen, Jean~C. Ragusa, Jim~E. Morel, Marvin~L. Adams, and Tan
  Bui-Thanh.
\newblock {Accelerating PDE-constrained Inverse Solutions with Deep Learning
  and Reduced Order Models}.
\newblock arXiv 1912.08864, 2019.

\bibitem{Shor12}
Naum~Zuselevich Shor.
\newblock {\em {Minimization Methods for Non-Differentiable Functions}},
  volume~3.
\newblock Springer Science \& Business Media, 2012.

\bibitem{SummWarrMoraLyge15}
Tyler Summers, Joseph Warrington, Manfred Morari, and John Lygeros.
\newblock Stochastic optimal power flow based on conditional value at risk and
  distributional robustness.
\newblock {\em International Journal of Electrical Power \& Energy Systems},
  72:116--125, 2015.

\bibitem{TomaRose14}
Vikrant~Singh Tomar and Richard~C Rose.
\newblock Manifold regularized deep neural networks.
\newblock In {\em Fifteenth Annual Conference of the International Speech
  Communication Association}, 2014.

\bibitem{VrakMathAnde14}
Maria Vrakopoulou, Johanna~L. Mathieu, and G{\"o}ran Andersson.
\newblock Stochastic optimal power flow with uncertain reserves from demand
  response.
\newblock In {\em 2014 47th Hawaii International Conference on System
  Sciences}, pages 2353--2362. IEEE, 2014.

\bibitem{Wils63}
R.~B. Wilson.
\newblock {\em {A Simplicial Algorithm for Concave Programming}}.
\newblock Ph.D. Thesis, Graduate School of Business Administration, Harvard
  University, Cambridge, MA, USA, 1963.

\bibitem{WoodWollSheb13}
Allen~J. Wood, Bruce~F Wollenberg, and Gerald~B Shebl{\'e}.
\newblock {\em Power generation, operation, and control}.
\newblock John Wiley \& Sons, 2013.

\bibitem{ZhuZabaKoutPerd19}
Yinhao Zhu, Nicholas Zabaras, Phaedon-Stelios Koutsourelakis, and Paris
  Perdikaris.
\newblock Physics-constrained deep learning for high-dimensional surrogate
  modeling and uncertainty quantification without labeled data.
\newblock {\em Journal of Computational Physics}, 394:56--81, 2019.

\end{thebibliography}

\end{document}